\documentclass[10.5pt,letter]{article}

\usepackage{booktabs, multirow}
\usepackage{fancyhdr}
\usepackage{tikz}
\usetikzlibrary{arrows.meta}
\usetikzlibrary{positioning}
\usepackage{bbm}
\usepackage{natbib}
\usepackage{amstext}
\usepackage{amsthm}
\usepackage{epsfig}
\usepackage{pdfpages}
\usepackage{amsmath,amssymb,amsthm}
\usepackage{graphicx}
\usepackage{calc}
\usepackage{dsfont}
\usepackage{rotating}
\usepackage{multirow}
\usepackage{amsmath}
\usepackage{lscape}
\usepackage{tabularx}
\usepackage{bm}
\usepackage{color}
\usepackage{enumerate,array}
\usepackage{verbatim}
\usepackage{stackengine}
\usepackage{enumitem}
\usepackage{mathrsfs}
\usepackage{titling}

\usepackage{hyperref}
\definecolor{cornellred}{rgb}{0.7, 0.11, 0.11}
\hypersetup{    
	colorlinks=true,
	linkcolor=blue,
	citecolor=blue,
	urlcolor=black
}

\renewcommand{\baselinestretch}{1.5}

\numberwithin{equation}{section}

\allowdisplaybreaks

\newtheorem{theorem}{Theorem}[section]
\newtheorem{lemma}[theorem]{Lemma}

\theoremstyle{definition}

\newtheorem{remark}[theorem]{Remark}

\newcommand{\R}{\mathbb{R}}

\newcommand{\N}{\mathbb{N}}

\newcommand{\E}{\mathbb{E}}
\renewcommand{\P}{\mathbb{P}}
\newcommand{\1}{\mathbbm{1}}

\newtheorem{innercustomgeneric}{\customgenericname}
\providecommand{\customgenericname}{}
\newcommand{\newcustomtheorem}[2]{%
  \newenvironment{#1}[1]
  {%
   \renewcommand\customgenericname{#2}%
   \renewcommand\theinnercustomgeneric{##1}%
   \innercustomgeneric
  }
  {\endinnercustomgeneric}
}
\newcustomtheorem{customthm}{Assumption}
\newcustomtheorem{customcor}{Corollary}

\renewcommand{\baselinestretch}{1.4}

\makeatletter
\def\singlespace{\deltaf\baselinestretch{1}\@normalsize}

\title{Sharp adaptive nonparametric testing for constant volatility }

\author{\textsc{Johannes Brutsche} and \textsc{Lukas Riepl} \\ Albert-Ludwigs-Universität Freiburg}
\date{}

\begin{document}

\maketitle

\vspace{-1cm}

\begin{abstract}
Based on discrete observations, we develop a test to infer if the volatility function $\sigma(\cdot)$ within the nonparametric Gaussian white noise model $dZ_t = \sigma(t)dW_t$ is constant. We propose a multiscale test based on localized squared increments, where the unknown volatility level under the composite null of constancy is estimated by quadratic variation. 
The testing procedure is shown to be minimax-optimal and attains the separation rate without knowledge of the smoothness parameters. For Hölder smoothness $0<\beta\leq 1$, we additionally determine the exact minimax constant and prove sharp adaptivity with respect to the Hölder radius. Our results identify the relative deviation $\sup_t|\sigma^2(t)/\|\sigma\|_{L^2}^2-1|$ as the natural scale-invariant measure of departure from the null hypothesis of constant volatility.
\end{abstract}

\vspace{0.5em}

\medskip 
\noindent {\bf Keywords:} 
Infill asymptotics, multiscale test, nonparametric volatility testing, sharp adaptivity.

\section{Introduction}

In this paper, we develop a (multiple) test to infer about constancy of the volatility function $\sigma$ in the nonparametric Gaussian white noise model 
\begin{align}\label{eq:Gaussian_white_noise}
    dZ_t = \sigma(t)dW_t, \qquad Z_0=z_0\in\R,\quad t\in [0,1],
\end{align}  
with driving noise $W$ being a Brownian motion and $\sigma:[0,1]\rightarrow (0,\infty)$ a deterministic, measurable and square-integrable function. More specifically, given high-frequency observations $Z_0,Z_{1/n},\dots, Z_{(n-1)/n}, Z_1$ of $(Z_t)_{t\in [0,1]}$, we develop a multiple test for the composite \textit{constancy testing problem}
\begin{align}\label{eq:testing_problem}
    H_0: \sigma \textrm{ is constant } \qquad\textrm{versus}\qquad H_1: \sigma \textrm{ is not constant}.
\end{align}

This test is shown to attain the minimax separation rate over Hölder classes in an adaptive way. Moreover, being a multiple test, it provides a data-driven localization of regions where deviations from constancy are indicated.


\begin{figure}[htbp]
    \centering
    \includegraphics[height=0.4\textwidth]{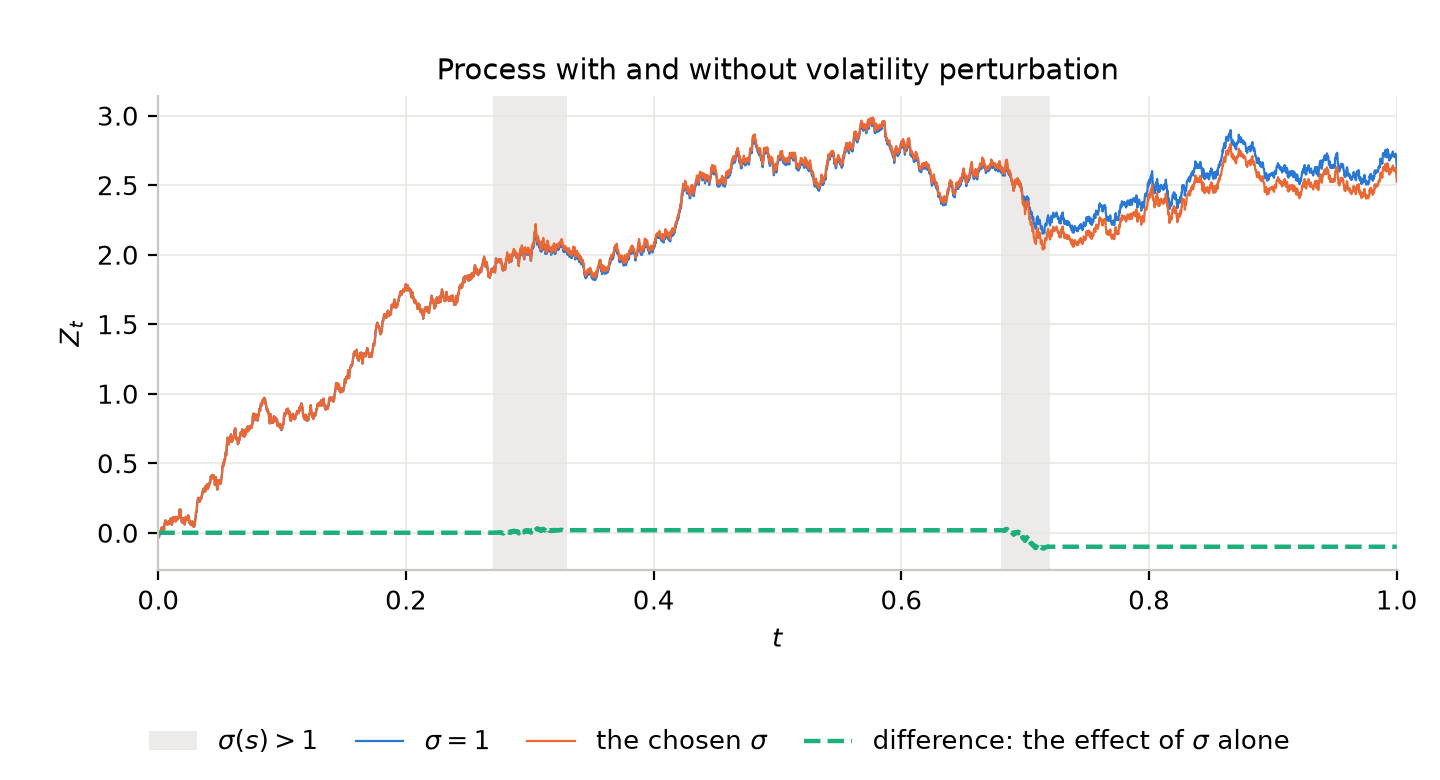} \vspace{-0.3cm}
    \caption{The sample path of a given Brownian motion (blue) together with the process $Z$ (orange) for some $\sigma=1+f$, where the support of $f$ is marked by the grey regions. In Section~\ref{section:simulation}, this example is analyzed using our  multiple testing procedure, which indeed detects both grey regions.}
    \label{fig:volatility_new}
\end{figure}

Generally speaking, the testing problem~\eqref{eq:testing_problem} is a special case of testing for the presence of a parametric function within a nonparametric model, which has been investigated in the context of regression models in~\cite{Mammen/Haerdle}. Motivated by financial econometrics, global goodness-of-fit tests for parametric classes were studied for volatility coefficients of stochastic differential equations in~\cite{Dette/Lieres/Wilkau}, \cite{Dette/Podolskij/Vetter} and~\cite{Dette/Podolskij}. In this article, we focus on the special parametric hypothesis of constancy for which we not only solve the testing problem~\eqref{eq:testing_problem}, but are additionally interested in the following statistical properties: First, we identify the minimax separation rate over Hölder classes for the problem~\eqref{eq:testing_problem} and develop a procedure that attains this rate. Secondly, we design this procedure to be  a multiple test statistic that is able to detect regions of deviations. The idea here is that for non-constant volatility $\sigma$, there exists at least one value $t\in [0,1]$ where $\sigma^2(t)$ differs from its squared $L^2$-norm $\|\sigma\|_{L^2}^2$. Given some regularity of $\sigma$, this possibly occurs on several unknown intervals and we develop our test in such a way that, in addition to the global test decision, it indicates where such intervals are located. This is illustrated in Figure~\ref{fig:volatility_new}.

The testing procedure we develop is based on localized squared increments and follows the so-called multiscale approach as described in~\cite{Dumbgen2001} in the context of testing qualitative hypotheses about the drift function in the Gaussian white noise model. Based on local log-likelihood ratio statistics, this approach naturally leads to a multiscale statistic suitable for testing $H_0: \sigma=c$ against $H_1: \sigma\neq c$ for a fixed constant $c$. To end up with a test statistic for the testing problem~\eqref{eq:testing_problem}, we propose to use the same statistic, but replace the unknown $c$ with the quadratic variation estimator. Incorporating this parametric estimator in the multiscale procedure is the main extension to existing multiscale tests in the literature. It will turn out that the resulting statistic is again tight and distribution free under $H_0$, where it is crucial that $\sigma=c$ can be estimated at the classical parametric rate by an estimator solely built upon squared increments of the observations. For sufficiently smooth volatility functions, our testing problem is also related to the Gaussian white noise framework of~\cite{Dumbgen2001} through asymptotic equivalence results such as those of~\cite{Grama/Nussbaum}. However, our direct approach covers the full range of Hölder exponents $\beta>0$, provides exact finite-sample calibration under the composite null hypothesis, and accounts for the non-Gaussian tails of squared increments across all bandwidths.

For the identification of the minimax separation rate, one has to quantitatively measure a deviation from the null hypothesis of constancy in~\eqref{eq:testing_problem}. Here, our study suggests working with the difference of the ratio of local volatility and its squared $L^2$-norm with one, i.e. measuring a deviation from constancy by $\sup_{t}|\sigma(t)^2 /\|\sigma\|_{L^2}^2 - 1|$. Like the hypothesis of constancy, this distance is scale-invariant and therefore naturally adapted to the problem, which is not the case for measuring deviations on absolute scale. For this deviation, we derive a lower bound over Hölder classes and finally establish that our testing procedure attains the optimal rate. By the special structure of being a multiscale test, we automatically adapt to the unknown smoothness parameters. Moreover, for Hölder smoothness $0<\beta\leq 1$, we additionally are optimal including the exact minimax constant and prove sharp adaptivity with respect to the Hölder radius.
While the alternatives we consider in this article consist of local but smooth departures from constancy, it should be noted that there exists a related strand of literature concerned with the location of structural breaks. Such change-point phenomena are for example discussed in~\cite{Bibinger/Jirak/Vetter}, \cite{Bibinger/Madensoy} and~\cite{Kokoszka}, see also references therein.

The article is organized as follows: In Section~\ref{Section:calibration} we develop our multiple test for constancy and provide details on the calibration of the multiscale test statistic (see Theorem~\ref{thm:tightness}). In Section~\ref{Section:power_properties}, we formally introduce the deviation measure from constancy and provide a detailed discussion. Following this, we present the minimax results for the particular choice of Hölder classes $\mathcal{H}(\beta,L)$, where the test is shown to be minimax optimal including the exact constant for all $0<\beta\leq 1$ and rate-optimal for every regime $\beta>1$. Moreover, the test is shown to be (sharp) adaptive in the unknown smoothness parameters $\beta$ and $L$ (see Theorems~\ref{thm:lower_bound}, \ref{thm:upper_bound} and~\ref{thm:adaptivity}). Section~\ref{section:simulation} contains a numerical illustration. A route of proof highlighting our main technical contributions is given after the respective results, while detailed proofs are provided in Section~\ref{Section:proofs}.

\section{The multiscale test for constant volatility}\label{Section:calibration}

In this section, we will develop a multiple test for the testing problem~\eqref{eq:testing_problem} based on the following vector of discrete observations
\[ \mathbf{Z}_n = \left( Z_0, Z_{1/n}, \dots, Z_{(n-1)/n}, Z_1\right).\]
Here and subsequently, the subscript $\sigma$ in the probability $\P_\sigma$ indicates that $Z$ is given as in~\eqref{eq:Gaussian_white_noise} with the corresponding volatility $\sigma$.
For the construction of the test statistic, suppose for a moment that we would like to test the simple null $\sigma=c$ against the simple alternative $\sigma=c+f$ for a fixed constant $c>0$ and a given deterministic, measurable and $L^2$-integrable function $f:[0,1]\rightarrow (-c,\infty)$. According to the Neyman-Pearson lemma, an optimal test is given by the log-likelihood ratio statistic $\log((\textrm{d}\P_{c+f}/\textrm{d}\P_c)(\mathbf{Z}_n))$. As $(Z_0,Z_{1/n},\dots, Z_1)$ is a Markov chain under $\P_\sigma$ where the conditional distribution 
\[ Z_{k/n}|Z_{(k-1)/n}\sim\mathcal{N}\left(Z_{(k-1)/n}, \int_{(k-1)/n}^{k/n}\sigma(s)^2ds\right) \]
is a Gaussian distribution, this log-likelihood ratio is given in our model as
\begin{align}\label{eq:likelihood_ratio}
    \sum_{k=1}^n \left(-\frac12 \log\left(\frac{n}{c^2}\int_{(k-1)/n}^{k/n}(c+f(s))^2ds\right) + \frac{n(Z_{k/n}-Z_{(k-1)/n})^2}{2c^2} \left[1-\frac{c^2}{n\int_{(k-1)/n}^{k/n}(c+f(s))^2ds}\right] \right). 
\end{align}
The logarithmic term does not depend on the data $\mathbf{Z}_n$ and is therefore omitted as we standardize the test statistic in one of the following steps by its mean. Assuming additionally that the support of $f$ is part of an interval $[t-h,t+h]$, the factor in square brackets is zero if $[(k-1)/n,k/n]\cap [t-h,t+h]=\emptyset$, meaning that this term acts like a localizer. For $h>0$ and a kernel $\psi:\R\rightarrow \R$ with support in $[-1,1]$ we define
\begin{align}\label{eq:sigma_n}
    \psi_{t,h}(s) := \psi\left(\frac{s-t}{h}\right), \quad \lambda_k(t,h) := \sqrt{n}\int_{(k-1)/n}^{k/n} \psi_{t,h}(s) ds \quad\textrm{ and }\quad\Lambda_n(t,h) := \sqrt{\sum_{k=1}^n \lambda_k(t,h)^2}.
\end{align}
Note that the dependence of $\lambda_k(t,h)$ and $\Lambda_n(t,h)$ on the kernel $\psi$ is suppressed in the notation. Moreover, in case $|\psi(x)|\leq 1$ for all $x\in\R$, we always have $\Lambda_n(t,h)\leq 1$. By the above mentioned heuristics we now substitute the local log-likelihood ratio statistic in~\eqref{eq:likelihood_ratio} by the local statistic
\begin{align}\label{eq:local_statistic_unnormalized}
    \sum_{k=1}^n \lambda_k(t,h) n\left(Z_{k/n}- Z_{(k-1)/n}\right)^2.
\end{align} 
Under $\P_c$, $(n(Z_{k/n}-Z_{(k-1)/n})^2/c^2)_{k=1,\dots,n}$ is a sequence of independent, identically $\chi_1^2$-distributed random variables. Thus, standardizing the expression in~\eqref{eq:local_statistic_unnormalized} by mean and variance provides the local test statistic
\begin{align}\label{eq:local_statistic_c}
    \hat{\Psi}_{n,t,h}^c(\mathbf{Z}_n) := \frac{\frac{1}{\sqrt{2}}\sum_{k=1}^n \lambda_k(t,h)\left( n(Z_{k/n}-Z_{(k-1)/n})^2 - c^2\right)}{c^2 \Lambda_n(t,h)}.
\end{align}
Based on these local test statistics we could now build a test for the nonparametric testing problem $H_0:\sigma=c$ versus $H_1:\sigma\neq c$. In order to move on to the composite null of constancy given in~\eqref{eq:testing_problem}, we substitute the unknown $c^2$ by the quadratic variation estimator
\begin{align}\label{eq:c_hat}
    \hat{c}_n^2 := \sum_{k=1}^n \left( Z_{k/n}-Z_{(k-1)/n}\right)^2.
\end{align}
This estimator is consistent, i.e.  $\hat{c}_n^2\rightarrow_{\P_\sigma} \int_0^1 \sigma(s)^2ds$,
by the Fisk approximation of quadratic variation (see Proposition~$18.17$ in~\cite{Kallenberg}) and the limit equals $c^2$ for $\sigma=c$ being constant. Moreover, an application of Markov's inequality reveals that it actually attains the classical parametric $1/\sqrt{n}$-convergence rate, see the proof of the subsequent Theorem~\ref{thm:tightness} in which this rate plays a key role.
The next result shows that the local statistics $\hat{\Psi}_{n,t,h}^{\hat{c}_n}(\mathbf{Z}_n)$ can be suitably combined for all $(t,h)$ in 
\begin{align*} 
    \mathcal{T}_n(K) := \left\{(t,h)\in\R^2: h\in [K/n,1/2) \textrm{ and } t\in [h,1-h] \right\}
\end{align*}
 to construct the desired multiple test. Note that for each $n$ and a continuous kernel $\psi$, the process $(t,h)\mapsto \hat{\Psi}_{n,t,h}^{\hat{c}_n}(\mathbf{Z}_n)$ is continuous.

\begin{theorem}[Calibration of the test statistic]\label{thm:tightness}
Let $\psi:\R\rightarrow \R_{\geq 0}$ be a function with $1=\psi(0)\geq\psi(x)$ for $x\neq 0$ and support in $[-1,1]$ that is of bounded variation with $|\psi(x)-\psi(y)|\leq |x-y|^\vartheta$ for some $0<\vartheta\leq 1$.  Then, there exists a constant $K_\psi$ such that $\Lambda_n(t,h)^2\geq 1/n$ for all $(t,h)\in\mathcal{T}_n(K_\psi)$. Let $\hat{\Psi}_{n,t,h}^{\hat{c}_n}(\mathbf{Z}_n)$ be given as in~\eqref{eq:local_statistic_c} for this kernel $\psi$ and $\hat{c}_n^2$ in~\eqref{eq:c_hat}. Then, for any  constant $c>0$, the family
\begin{align*}
  \left( \sup_{(t,h)\in\mathcal{T}_n(K_\psi)} \left(\big|\hat{\Psi}_{n,t,h}^{\hat{c}_n}(\mathbf{Z}_n)\big| - G_n\left(2\log(1/\Lambda_n(t,h)), \Lambda_n(t,h)\right)\right)\right)_{n\in\N}  
\end{align*}
is tight under $\P_c$, where
\begin{align*}
    G_n(\eta,\delta) := \sqrt{2\eta}+\frac{2\eta}{\max\{1,\sqrt{n}\delta\}}.
\end{align*}
\end{theorem}

The assumption on the Hölder continuity of the kernel $\psi$ is needed for the bound $\Lambda_n(t,h)^2\geq 1/n$ and allows us to control the minimum number of time points $k/n$ within the support of the rescaled kernel $\psi_{t,h}$, see Lemma~\ref{lemma:sigmageq} for a detailed derivation.
The calibration of the test statistic relies crucially on distinguishing between the bias generated by taking the supremum over $(t,h)$ that necessitates the additive correction $G_n$, and the contribution of the parametric estimator $\hat{c}_n$ of the unknown $c$. The key inequality here is
\begin{align}\label{eq:tightnessleftright}
  & \sup_{(t,h)\in\mathcal{T}_n(K_\psi)} \left(\big|\hat{\Psi}_{n,t,h}^{\hat{c}_n}(\mathbf{Z}_n)\big| - G_n\left(2\log(1/\Lambda_n(t,h)), \Lambda_n(t,h)\right)\right) \begin{cases} \leq \Psi_n^1(\mathbf{Z}_n) + \Psi_n^2(\mathbf{Z}_n), \\ \geq  \Psi_n^2(\mathbf{Z}_n)- \Psi_n^1(\mathbf{Z}_n), \end{cases}
\end{align}
where
\begin{align*}
    \Psi_n^1(\mathbf{Z}_n) &= \sup_{(t,h)\in\mathcal{T}_n(K_\psi)} \left(\big|\hat{\Psi}_{n,t,h}^{c} (\mathbf{Z}_n)-\hat{\Psi}_{n,t,h}^{\hat{c}_n}(\mathbf{Z}_n) \big| \right),\\
    \Psi_n^2(\mathbf{Z}_n) &= \sup_{(t,h)\in\mathcal{T}_n(K_\psi)} \left(\big|\hat{\Psi}_{n,t,h}^{c} (\mathbf{Z}_n)\big| - G_n\left(2\log(1/\Lambda_n(t,h)), \Lambda_n(t,h)\right)\right).
\end{align*}
The term $\Psi_n^1$ forms a difference to existing work and treats the effect of including the parametric estimator $\hat{c}_n$ into the multiscale statistic. The proof of tightness of this term relies on the parametric $1/\sqrt{n}$-convergence rate of $\hat{c}_n^2$ towards $c^2$ that allows for an explicit upper bound of $\Psi_n^1$ not involving a supremum any longer. On the other hand, $\Psi_n^2$ equals our test statistic with known value $c$ of the volatility level and can be analyzed under $\P_c$ using arguments already present in other multiscale testing procedures in the literature, e.g.~\cite{Dumbgen_Annals}, \cite{Dumbgen_report} testing shape constraints of densities, \cite{Rohde2009} covering higher dimensions or~\cite{Brutsche/Rohde} with a test on the drift of an ergodic diffusion model. 
The function $G_n(\eta,\delta)$ is identified for our model as an inverted Bernstein tail which is due to the fact that our test statistic is built upon the squared increments $n(Z_{k/n}-Z_{(k-1)/n})^2/c^2$ that follow a $\chi_1^2$-distribution under $\P_c$. The details are provided in Subsection~\ref{Subsection:proof_calibration}.\smallskip

We now return to defining a test for the testing problem~\eqref{eq:testing_problem} on the basis of the statistic in Theorem~\ref{thm:tightness}. The vector $(n(Z_{k/n}-Z_{(k-1)/n})^2/c^2)_{k=1,\dots, n}$ consists of independently $\chi_1^2$-distributed random variables and both $\hat{c}_n^2/c^2$ and $\hat{\Psi}_{n,t,h}^{c}(\mathbf{Z}_n)$ are continuous functions of this vector, where $t,h$ only appear in deterministic factors. Therefore, it is easy to see that the process
\[ \left( \hat{\Psi}_{n,t,h}^{\hat{c}_n}(\mathbf{Z}_n)\right)_{(t,h)\in\mathcal{T}_n(K_\psi)}\]
has the same distribution under $\P_1$ and $\P_c$ for every constant $c>0$. For $\alpha\in (0,1)$, we define
\begin{align*}
    \kappa_n^\alpha := \inf\left\{ r\in\R: \P_1\left( \sup_{(t,h)\in\mathcal{T}_n(K_\psi)}\left( \big|\hat{\Psi}_{n,t,h}^{\hat{c}_n}(\mathbf{Z}_n)\big| - G_n\left(2\log(1/\Lambda_n(t,h)), \Lambda_n(t,h)\right)\right) \leq r\right)\geq 1-\alpha \right\},
\end{align*}
for which Theorem~\ref{thm:tightness} entails $\limsup_{n\to\infty}|\kappa_n^\alpha|<\infty$. Finally, we introduce the \textit{test for constant volatility}
\begin{align}\label{eq:test}
    T_n(\mathbf{Z}_n) := \1_{\left\{\sup_{(t,h)\in\mathcal{T}_n(K_\psi)}\left( \big|\hat{\Psi}_{n,t,h}^{\hat{c}_n}(\mathbf{Z}_n)\big| - G_n\left(2\log(1/\Lambda_n(t,h)), \Lambda_n(t,h)\right)\right) >\kappa_n^\alpha\right\}}
\end{align}
and as $(\hat{\Psi}_{n,t,h}^{\hat{c}_n}(\mathbf{Z}_n))_{(t,h)\in\mathcal{T}_n(K_{\psi})}$ is distribution free under the null hypothesis, we directly obtain that the test is uniformly (over $H_0$) of level $\alpha$, i.e. that for every $n\in\N$,
\[ \sup_{c\textrm{ constant}}\P_c\left( T_n(\mathbf{Z}_n) =1\right) \leq \alpha.\]
We highlight at this point that the above statement is non-asymptotic. The power properties of the test for constant volatility are studied in the subsequent Section~\ref{Section:power_properties} for the particular example of Hölder smoothness classes. It also contains the discussion on how to appropriately measure deviations from the null hypothesis of constancy in~\eqref{eq:testing_problem}.

\begin{remark}[Simultaneous detection of regions of deviation]\label{remark:D_alpha}
The test $T_n(\mathbf{Z}_n)$ rejects $H_0$ if, and only if, the random family
\[ \mathcal{D}_n^\alpha := \left\{ (t,h)\in\mathcal{T}_n(K_\psi): \big|\hat{\Psi}_{n,t,h}^{\hat{c}_n}(\mathbf{Z}_n)\big| > \kappa_n^\alpha + G_n\left( 2\log(1/\Lambda_n(t,h)), \Lambda_n(t,h)\right)\right\}\]
is non-empty. In particular, $\sup_{c\textrm{ constant}}\P_c(\mathcal{D}_n^\alpha\neq\emptyset)\leq\alpha$. The intervals $[t-h,t+h]$ with $(t,h)\in\mathcal{D}_n^\alpha$ indicate where the data deviate from constancy and thus provide a data-driven localization of potential deviations of $\sigma^2$ from $\|\sigma\|_{L^2}^2$, see Section~\ref{section:simulation} for an illustration.
\end{remark}

\section{Minimax optimality and sharp adaptivity}\label{Section:power_properties}

In this section we establish minimax results for the test for constant volatility given in~\eqref{eq:test}. This is exemplarily done for Hölder smoothness classes which are introduced below. As the test is built as a multiple test taking into account all possible bandwidths, we also derive that it is automatically adaptive in the unknown smoothness parameters.

\subsection{Measuring deviations from constancy}
In order to determine minimax rates, we need to specify a distance between a given volatility function $\sigma$ and the class of constant functions, i.e. those lying in the null hypothesis of our testing problem. Heuristically, an absolute perturbation of fixed size is harder to detect the higher the baseline volatility. Deviations from constancy should therefore be measured relative to the global volatility level rather than on an absolute scale. We define a distance as
\begin{align*}
    d_J(\sigma,H_0) := \sup_{t\in J}\left| \frac{\sigma(t)^2}{\|\sigma\|_{L^2}^2} -1 \right|
\end{align*}
for compact $J\subset (0,1)$ to avoid boundary effects. This means that local deviations from the null at minimax detection scale are not truncated by the endpoints of the observation domain, which allows for the derivation of optimal constants. Note that $d_J(\sigma,H_0)=0$ in case $\sigma=c$ is constant and that $d_J(\sigma,H_0)>0$ for $\sigma$ being non-constant on $J$. Moreover, our particular choice of $d_J$ reflects the fact that constancy of $\sigma$ is a scale-invariant property, and testing for constancy should be of the same level of difficulty for data emerging from the process $Z$ or $aZ$ for any $a>0$.

\subsection{Optimal power properties}
For an interval $I\subset\R$ and $\beta,L>0$ denote by $\mathcal{H}(\beta,L)$ the H\"older class of functions, that is all functions $f:I\rightarrow\R$ for which the first $\lfloor\beta\rfloor$ derivatives exist and
\[ \left| f^{(\lfloor \beta\rfloor)}(x) -f^{(\lfloor \beta\rfloor)}(y)\right| \leq L|x-y|^{\beta-\lfloor\beta\rfloor}, \]
where $f^{(k)}$ denotes the $k$th derivative of $f$ and $\lfloor\beta\rfloor$ the maximal integer strictly smaller than $\beta$. The domain $I$ will always be clear from the context. Our subsequent power results will specify those volatility functions $\sigma$ satisfying $\sigma^2/\|\sigma\|_{L^2}^2\in\mathcal{H}(\beta,L)$ that are close enough to the null hypothesis of constancy in the distance $d_J$ such that they can be detected by the test $T_n$ given in~\eqref{eq:test}. To this aim, define the rate and constant
\begin{align}\label{eq:optimal_constant}
    \rho_n = \rho_n(\beta):=\left(\frac{\log n}{n}\right)^{\frac{\beta}{2\beta+1}} \qquad\textrm{ and }\qquad c_* = c_*(\beta,L) :=\left( \frac{4L^{1/\beta}}{(2\beta+1)\|\psi_\beta\|_{L^2}^2}\right)^{\frac{\beta}{2\beta+1}},
\end{align} 
where $\psi_\beta$ is the so-called \textit{optimal recovery kernel} which by definition is the solution to the optimization problem
\[ \textrm{ Minimize }\|\psi\|_{L^2} \textrm{ over all } \psi\in\mathcal{H}(\beta,1)\textrm{ with } \psi(0)\geq 1.\]
It is known that $\psi_\beta$ is an even function with compact support and satisfies $\psi_\beta(0)=1>|\psi_\beta(x)|$ for $x\neq 0$, see for example~\cite{Leonov} or~\cite{Donoho}. For $0<\beta\leq 1$, it has the explicit representation
\begin{align}\label{eq:optimal_recovery_kernel}
    \psi_\beta(x) = \1_{\{|x|\leq 1\}}\left(1-|x|^\beta\right).
\end{align} 

We start the minimax theory by providing a lower bound. The following result establishes that for every test of level $\alpha$ for the null hypothesis of constancy there exist volatility coefficients $\sigma$ in the alternative that deviate $(1-\epsilon_n)c\rho_n$ from being constant in the distance $d_J$ which cannot be detected with probability $1-\alpha-o(1)$ or larger. Note that this result even holds true under knowledge of both smoothness parameters $\beta$ and $L$.

\begin{theorem}[Minimax lower bound]\label{thm:lower_bound}
   Let $\beta,L>0$, $\alpha\in (0,1)$ and $(\phi_n)_{n\in\N}$ a sequence of tests satisfying $\sup_{\sigma\in H_0} \E_\sigma[\phi_n(\mathbf{Z}_n)]\leq \alpha$. Then for arbitrary numbers $\epsilon_n>0$ with $\lim_{n\to\infty}\epsilon_n=0$ and $\lim_{n\to\infty}\epsilon_n^2 \log(n)=\infty$ there exists a constant $c=c(\beta,L)$ such that 
   \[ \limsup_{n\to\infty} \inf_{ \substack{\sigma^2/\|\sigma\|_{L^2}^2\in\mathcal{H}(\beta,L):\\ d_J(\sigma, H_0)\geq (1-\epsilon_n)c\rho_n}} \E_\sigma[\phi_n(\mathbf{Z}_n)]\leq \alpha \]   
   for any fixed compact interval $J=[J_1,J_2]\subset (0,1)$ with $J_1<J_2$.
In case $\beta\in (0,1]$ the result holds true for $c=c_*$.
\end{theorem}

The proof relies crucially on the construction of hypotheses $\sigma_{n,j}$ for $j=1,\dots, N_n$ and $N_n$ appropriate in such a way that
\begin{itemize}
    \item[(1)] $\sigma_{n,j}^2/\|\sigma_{n,j}\|_{L^2}^2\in\mathcal{H}(\beta,L)$ and
    \item[(2)] $d_J(\sigma_{n,j},H_0) \geq (1-\epsilon_n)c\rho_n$ with $c=c_*$ for $\beta\leq 1$.
\end{itemize}
Here, item (2) is particularly challenging, since the deviation measure $d_J$ compares the local quantity $\sigma_{n,j}(t)$ with the global quantity $\|\sigma_{n,j}\|_{L^2}$. Although we construct $\sigma_{n,j}$ in the common way of being local deviations from some baseline $\sigma_0$, the scaling required for these to fulfill (2) is not a priori clear. In view of (1), a natural choice for $\sigma_{n,j}$ requires the form
\begin{align}\label{eq:hypothesis_squared}
    \sigma_{n,j}^2(t) = 1+ L(1-\epsilon_n) h_n^\beta \psi\left(\frac{t-t_{n,j}}{h_n}\right)
\end{align} 
for some kernel $\psi$, suitable bandwidth $h_n$ and location parameters $t_{n,j}$. Indeed, with such a choice we can obtain the rate $\rho_n$ and we employ it for $\beta>1$. Deriving the optimal constant $c_*$, however, requires another choice. Here, we introduce a parameter $b>0$ and define a bandwidth
\begin{align}\label{eq:h_nb}
    h_n^b := \left(\frac{c_*}{L}\right)^{\frac{1}{\beta}} \left( \frac{\log(n)}{bn}\right)^{\frac{1}{2\beta+1}}.
\end{align} 
For location parameters $t_{n,j}$ that will be specified later in the detailed proof, we introduce the hypotheses
\begin{align}\label{eq:sigma_j}
 \sigma_{n,j}^b(t) = 1+\frac{L}{2}(1-\epsilon_n)(h_n^b)^\beta \psi_\beta\left(\frac{t-t_{n,j}^b}{h_n^b}\right).
\end{align} 
As compared to~\eqref{eq:hypothesis_squared}, the important difference is that we now specify $\sigma_{n,j}$ in the given way and not $\sigma_{n,j}^2$.
The following two lemmas establish that there exists a sequence $(b_n)_{n\in\N}$ such that items (1) and (2) are satisfied. Moreover, the second result also contains a statement about the speed of convergence of the solutions $b_n$ towards one, which is essential to derive the exact constant in the proof of Theorem~\ref{thm:lower_bound}.

\begin{lemma}\label{lemma:hypotheses_Hölder}
    Let $\beta\leq 1$, $(b_n)_{n\in\N}$ a sequence bounded away from zero and $\sigma_{n,j}^{b_n}$ be given as in~\eqref{eq:sigma_j}. Then there exists $n_0\in\N$ such that $(\sigma_{n,j}^{b_n})^2/\|\sigma_{n,j}^{b_n}\|_{L^2}^2 \in\mathcal{H}(\beta,L)$ for all $n\geq n_0$.
\end{lemma}

\begin{lemma}\label{lemma:hypotheses_solution}
    Let $\beta\leq 1$ and $\sigma_{n,j}^b$ be given as in~\eqref{eq:sigma_j}. Then there exists $n_0$ and $b=b_n>0$ such that $\sigma_{n,j}^{b_n}$ solves (2) for $n\geq n_0$. Moreover, we have $\log(n)|b_n-1|\rightarrow 0$.
\end{lemma}

The strategy to prove Lemma~\ref{lemma:hypotheses_solution} is to define $b_n$ as a solution of the equality in
\begin{align}\label{eq_proof:hypotheses_solution1}
    d_J(\sigma_{n,j}^b,H_0) \geq \left|\frac{\sigma_{n,j}^b(t_{n,j})^2}{\|\sigma_{n,j}^b\|_{L^2}^2} -1\right| = (1-\epsilon_n)c_*\rho_n.
\end{align}  
Both the existence of such $b_n$ and the bound on $|b_n-1|$ can be derived via the implicit function theorem. All details for this, the proof of Lemma~\ref{lemma:hypotheses_Hölder} and the proof of Theorem~\ref{thm:lower_bound} are provided in Subsection~\ref{Subsection:proof_lower_bound}. \smallskip

The next result shows that our test $T_n$ given in~\eqref{eq:test} achieves the corresponding rate to Theorem~\ref{thm:lower_bound} and in case $\beta\leq 1$ also the optimal constant.

\begin{theorem}[Minimax upper bound]\label{thm:upper_bound}
   Let $\beta,L>0$ and $\psi$ a kernel satisfying the assumptions of Theorem~\ref{thm:tightness}. Then for arbitrary numbers $\epsilon_n>0$ with $\lim_{n\to\infty}\epsilon_n=0$ and $\lim_{n\to\infty} \epsilon_n^2\log(n)=\infty$ there exists a constant $c=c(\beta,L,\psi)$ such that for the test $T_n$ given in~\eqref{eq:test},
   \[ \lim_{n\to\infty} \inf_{ \substack{\sigma^2/\|\sigma\|_{L^2}^2\in\mathcal{H}(\beta,L):\\ d_J(\sigma, H_0)\geq (1+\epsilon_n)c\rho_n}}\P_\sigma\left( T_n(\mathbf{Z}_n) =1\right)=1\]
   for any fixed compact interval $J=[J_1,J_2]\subset (0,1)$ with $J_1<J_2$. In case $\beta\in (0,1]$, we can choose $\psi=\psi_\beta$ and the result holds true for $c=c_*$.
\end{theorem}

From the last result we know that any $\sigma$ generating the data $\mathbf{Z}_n$ and deviating at least $(1+\epsilon_n)c\rho_n$ in the distance $d_J$ from the null hypothesis of constancy can be detected (i.e. the null is rejected) by the test $T_n$ given in~\eqref{eq:test} with probability close to one. In particular, this test $T_n$ attains the minimax optimal rate over the whole range of $\beta,L>0$ and additionally is even optimal in the constant for $\beta\leq 1$. The next result states that rate-adaptivity is attained even uniformly over $(\beta,L)\in [\beta_1,\beta_2]\times[L_1,L_2]$ and sharp adaptivity for fixed $\beta$ over $L\in [L_1,L_2]$.

\begin{theorem}[Adaptivity]\label{thm:adaptivity}
    Let $\psi,\epsilon_n$ be specified as in Theorem~\ref{thm:upper_bound}. Then for any compact subset $[\beta_1,\beta_2]\times[L_1,L_2]\subset (0,\infty)^2$ the test $T_n$ in~\eqref{eq:test} is rate-adaptive in both parameters $\beta$ and $L$ in the sense that for $J=[J_1,J_2]\subset (0,1)$ with $J_1<J_2$,
    \[ \lim_{n\to\infty} \inf_{(\beta,L)\in [\beta_1,\beta_2]\times[L_1,L_2]}\inf_{ \substack{\sigma^2/\|\sigma\|_{L^2}^2\in\mathcal{H}(\beta,L):\\ d_J(\sigma, H_0)\geq (1+\epsilon_n)c\rho_n}}\P_\sigma\left( T_n(\mathbf{Z}_n) =1\right)=1.\]
    For $\beta\leq 1$  and $T_n$ with $\psi=\psi_\beta$, we have sharp adaptivity in the parameter $L\in[L_1,L_2]\subset (0,\infty)$ in the sense
    \[ \lim_{n\to\infty} \inf_{L\in [L_1,L_2]}\inf_{ \substack{\sigma^2/\|\sigma\|_{L^2}^2\in\mathcal{H}(\beta,L):\\ d_J(\sigma, H_0)\geq (1+\epsilon_n)c_*\rho_n}}\P_\sigma\left( T_n(\mathbf{Z}_n) =1\right)=1.\]
\end{theorem}

\begin{remark}[Choice of kernel]
For any kernel $\psi$ satisfying the assumptions of Theorem~\ref{thm:tightness}, the multiscale test~\eqref{eq:test} is rate-optimal and adaptive by Theorems~\ref{thm:upper_bound} and~\ref{thm:adaptivity}. The specific choice of the optimal recovery kernel $\psi_\beta$ is only required to attain the exact detection constant and requires the smoothness $0<\beta\leq 1$ to be known.
\end{remark}

\begin{remark}[Optimal constant]
For composite hypotheses testing about the drift function in Gaussian white noise models, the optimal constant was derived in~\cite{Dumbgen2001}. Here, $4$ in the constant in~\eqref{eq:optimal_constant} is replaced with $2$. Our additional factor $2$ is due to the fact that $\mathrm{Var}(Y)=2$ for $Y\sim\chi_1^2$ and our statistic being built upon the quantities $n(Z_{k/n}-Z_{(k-1)/n})^2/c^2$, which follow a $\chi_1^2$-distribution under $\P_c$, whereas the drift is estimated based on the increments $\sqrt{n}(Z_{k/n}-Z_{(k-1)/n})/c\sim\mathcal{N}(0,1)$.
\end{remark}

\begin{remark}[Possible extensions]
We briefly discuss two possible extensions of our framework. We regard both of them as natural directions for further research, while the present paper focuses on the setting in which exact finite-sample calibration and sharp adaptive minimax optimality can be established simultaneously.
\begin{itemize}
\item \textbf{Stochastic volatility.} For state-dependent volatility $\sigma(t,Z_t)$, our test retains its exact calibration under constant volatility in the driftless model. Under alternatives, however, the increments are no longer independent Gaussian variables with deterministic variances, and deviations may occur at random locations and persist for random durations. Their detectability thus depends on the occupation of relevant state-space regions. Together with the different temporal regularity of realized volatility, this requires a new formulation of the minimax testing problem.
\item \textbf{General parametric hypotheses.} For deterministic volatility models $H_0=\{\sigma^2=q_\theta:\theta\in\Theta\}$, standardized squared increments remain independent $\chi_1^2$-distributed variables when $\theta$ is known, so the basic multiscale construction carries over. Unknown parameters require uniform estimation bounds and potentially a different calibration. Moreover, unlike the hypothesis of constancy, where scale invariance motivates the distance $d_J$, general parametric families do not come with a natural separation distance, and different projections onto the null model may lead to different sharp constants.
\end{itemize}
\end{remark}

\section{Numerical illustration}\label{section:simulation}
\begin{figure}[!htbp]
        \centering
        \includegraphics[width=0.6\textwidth]{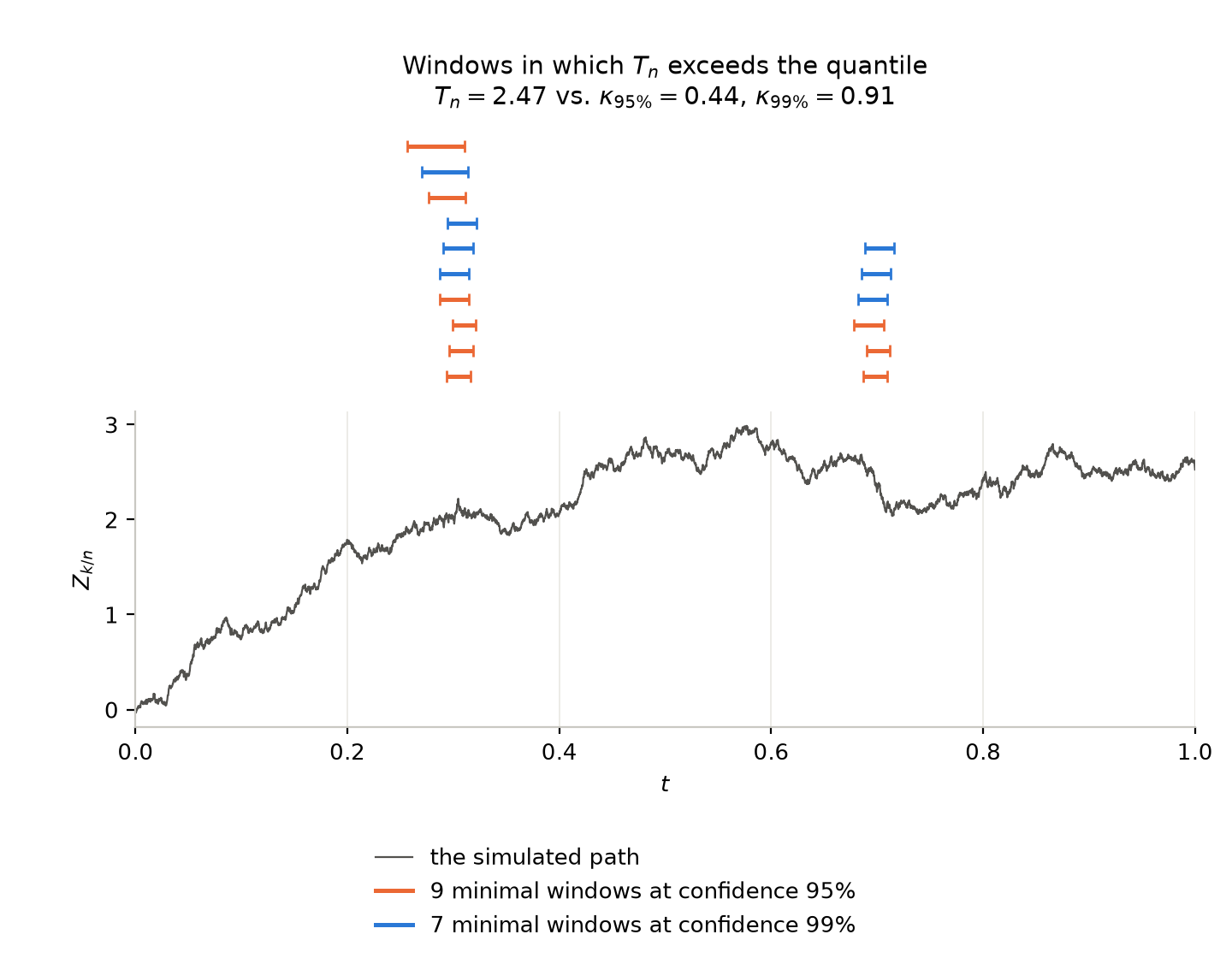}\vspace{-0.3cm}
        \caption{Visualization of the sets $\mathcal{D}^{0.05}_{5000}$ and  $\mathcal{D}^{0.01}_{5000}$}
        \label{fig:detected regions}
\end{figure}

\begin{figure}[!htbp]
        \centering
        \includegraphics[width=0.5\textwidth]{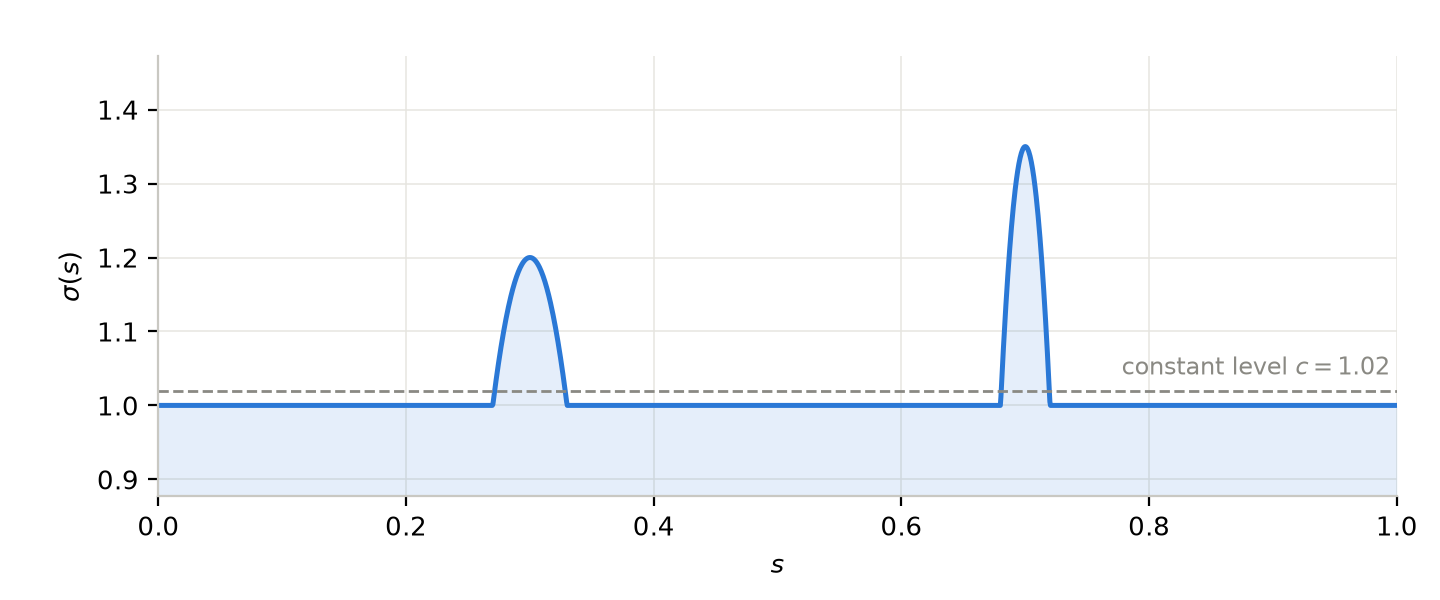}\vspace{-0.3cm}
        \caption{$\sigma(s)$ and its $L^2$-norm}
        \label{fig:sigma(s)}
\end{figure}

We illustrate our test for the example already introduced in Figure~\ref{fig:volatility_new}. The observations are based on the volatility function
$$\sigma(s)=1+0.2\max\left\{0,1-\left(\left(s-0.3\right)/0.03\right)^2 \right\}+0.35\max\left\{0,1-\left(\left(s-0.7\right)/0.02\right)^2\right\}$$
and are simulated for $n=5000$. The specific choice of volatility function is illustrated in Figure~\ref{fig:sigma(s)}.
The test statistic is computed as a maximum over  the discretized grid
\[ \mathcal{T} = \left\{ (t,h)\in [0,1]\times [0,1/2]: t\leq 1-h, h= 6\cdot(1.25)^k/5000, t=(1+0.25\cdot l)h, \textrm{ with } k,l\in\N_0   \right\}\]
and all integrals in test statistic are approximated as $\int_{(k-1)/n}^{k/n}\psi_{t,h}(s)ds \approx \frac{1}{n}\psi_{t,h}(\frac{k-1/2}{n})$ to keep the runtime minimal. As a kernel, we choose the optimal recovery kernel from~\eqref{eq:optimal_recovery_kernel} for $\beta=1$. The quantiles $\kappa_{5000}^{0.05}$ and $\kappa_{5000}^{0.01}$ are estimated using $9999$ Monte-Carlo simulations. 
Recall the regions of deviations $\mathcal{D}_n^\alpha$ given in Remark~\ref{remark:D_alpha}. Figure~\ref{fig:detected regions} shows that the test successfully identifies the two localized perturbations and visualizes the minimal intervals in $\mathcal{D}_{5000}^{0.05}$ and $\mathcal{D}_{5000}^{0.01}$. Here, an interval $I\in \mathcal{D}_n^\alpha$ is called minimal, if for any interval $I\neq J\in\mathcal{D}_n^\alpha$ we have $J\nsubseteq I$.

\vspace{-0.4cm}

\section{Proofs}\label{Section:proofs}

This section is organized as follows: First we introduce some additional notation and present two preliminary lemmas that are used throughout the proofs. In Subsection~\ref{Subsection:proof_calibration} we then present the proof of Theorem~\ref{thm:tightness} and afterwards the proofs of the minimax results are given, where Subsection~\ref{Subsection:proof_lower_bound} contains the proof of the lower bound (Theorem~\ref{thm:lower_bound}) and Subsection~\ref{Subsection:Proof_upper_bound} the corresponding upper bound and adaptivity result (Theorem~\ref{thm:upper_bound} and Theorem~\ref{thm:adaptivity}).\smallskip

\textbf{Notation:} Here and throughout, we use the Landau symbols $o(a_n)$ and $\mathcal O(a_n)$, where $b_n=o(a_n)$ means $|b_n/a_n|\to 0$ as $n\to\infty$ and $b_n=\mathcal{O}(a_n)$ that $\limsup_n |b_n/a_n| <\infty$. Their stochastic counterparts $o_{\P}(1)$ and $\mathcal O_{\P}(1) $ denote convergence to zero and boundedness in probability, respectively. For any function $f:\R\rightarrow\R$, we denote $\|f\|_{\sup}:= \sup_{t\in\R}|f(t)|$ and as in the previous sections, $\|f\|_{L^2}$ denotes the $L^2$-norm, i.e. $\|f\|_{L^2}^2= \int_\R f^2(x)dx$.  \smallskip

The following lemma is a slight modification and reformulation of a result established in \cite{Rohde_2008} within the proof of Theorem 3, page 1371--1372.

\begin{lemma}[see \cite{Rohde_2008}, p.~1372]\label{lem:rohde2}
Let $\beta,L>0$ and $(f_n)_{n\in\mathbb{N}}$ be a sequence of functions defined on $[-a,a]$ with $f_n\in\mathcal{H}(\beta,L)$ and
$|f_n(x^*)| = \sup_{t\in [-a,a]}|f_n(t)| = \rho_n \le 1$. Then there exists a constant $c=c(\beta,L)>0$ such that
\[ |f_n(x)| \ge (1/2)\rho_n \qquad \text{for} \qquad x\in [-a,a]\cap \left[x^*-c\rho_n^{1/\beta},\,x^*+c\rho_n^{1/\beta}\right]. \]
Moreover, considered as a function in $\beta$ and $L$, the constant
$c(\beta,L)$ is continuous in its parameters.
\end{lemma}

The next result establishes that $\Lambda_n(t,h)^2$ is bounded from below by a multiple of the bandwidth $h$ in case $h$ is not too small. 

\begin{lemma}\label{lemma:sigmageq}
Let $\psi$ be given as in Theorem~\ref{thm:tightness}. Then there exists a constant $K_\psi$ such that $\Lambda_n(t,h)^2\geq 1/n$ for all $(t,h)\in\mathcal{T}_n(K_\psi)$. Moreover, there exists a constant $C=C(\psi)$ such that
$$ \Lambda_n(t,h)^2 \geq Ch  \quad\textrm{ for } (t,h)\in\mathcal{T}_n(K_\psi).$$
\end{lemma}
\begin{proof}
Using that $\psi$ is nonnegative allows us to bound 
\begin{align*}
\Lambda_n(t,h)^2 & \geq \frac{n}{4}\sum_{k=1}^n  \left(\int_{(k-1)/n}^{k/n}\1_{\{ \psi_{t,h}(s)  \geq 1/2 \}}ds\right)^2 .
\end{align*}
As $\psi(0)=1$, we obtain from the Hölder continuity assumption $|\psi(x)-\psi(0)|\leq |x|^\vartheta$ that $\psi(x)\geq 1/2$ for $|x|\leq 2^{-1/\vartheta}$. In particular, the support of $\1_{\{  \psi_{t,h}(s) \geq 1/2 \}}$ contains $[t-2^{-1/\vartheta} h, t+2^{-1/\vartheta} h]$. The number of indices $k$ for which $[(k-1)/n,k/n]\subset [t-2^{-1/\vartheta} h, t+2^{-1/\vartheta} h]$ is at least $ \lfloor 2^{1-1/\vartheta} nh \rfloor -1 \geq  2^{1-1/\vartheta}n h-2 $ such that 
\[ \Lambda_n(t,h)^2 \geq \frac{1}{2n}\left( 2^{-1/\vartheta} nh - 1\right). \]
Setting $K_\psi:= 3\cdot 2^{1/\vartheta}$ then yields $\Lambda_n(t,h)^2\geq 1/n$ for $h\geq K_\psi/n$. Moreover, 
for these $h$ we have $1/n\leq h/K_\psi$ and obtain the second assertion as
\[ \Lambda_n(t,h)^2 \geq \frac12 \left( 2^{-1/\vartheta} h - \frac1n\right) \geq \frac12\left( 2^{-1/\vartheta} -\frac{1}{K_\psi}\right) h = \frac13 2^{-1/\vartheta} h.\]
\end{proof}

\subsection{Proofs of Section~\ref{Section:calibration} (Calibration)}\label{Subsection:proof_calibration}

The identification of the correction term $G_n$ in Theorem~\ref{thm:tightness} is based on the following result on the tail probabilities of weighted sums of independent centered $\chi_1^2$-distributed random variables. This tail bound was derived in~\cite{Massart} and is given below in the form that we will use in the proof of Theorem~\ref{thm:tightness} for the reader's convenience.

\begin{lemma}[cf. Lemma~$1$ in~\cite{Massart}]\label{lemma:Bernstein_chi_squared}
Let $X_1,\dots, X_n$ be independent and identically standard Gaussian random variables and $a_1,\dots, a_n\in\R$. Define $V_n := \sum_{k=1}^n |a_k|^2$ and $M_n:=\max_{k=1,\dots,n}|a_k|$.
Then for any $x\geq 0$,
\begin{align*}
    \P\left( \left| \sum_{k=1}^n a_k \big(X_k^2-1\big)\right| \geq 2\sqrt{V_n x} + 2M_n x\right) \leq 2\exp\left( -x\right).
\end{align*}
\end{lemma}
\begin{proof}
For $u\in(0,1/2)$ we have $\log\E[e^{u(X_1^2-1)}]\leq u^2/(1-2u)$, see Lemma~$1$ in~\cite{Massart}, and $\log\E[ e^{-u(X_1^2-1)}]=u-\frac12\log(1+2u)\le u^2$. Hence, for $S=\sum_k a_k(X_k^2-1)$ and $0<\lambda<1/(2M_n)$, both $\log\E[e^{\lambda S}]$ and $\log\E[e^{-\lambda S}]$ are bounded by $\lambda^2V_n/(1-2\lambda M_n)$. Choosing $\lambda=\sqrt{x}/\big(\sqrt{V_n}+2M_n\sqrt{x}\big)\in(0,1/(2M_n))$ in Chernoff's bound $\P(\pm S\geq t)\leq \exp(-\lambda t)\E e^{\pm\lambda S}$ with $t=2\sqrt{V_nx}+2M_nx$ yields $\lambda t-\lambda^2V_n/(1-2\lambda M_n)=x$, and the claim follows by a union bound.
\end{proof}

\begin{proof}[Proof of Theorem \ref{thm:tightness}]
The existence of $K=K_\psi$ with $\Lambda_n(t,h)^2\geq 1/n$ for $(t,h)\in\mathcal{T}_n(K)$ was already proven in Lemma~\ref{lemma:sigmageq}.
Following the route of proof after Theorem~\ref{thm:tightness}, it remains to prove that both quantities $\Psi_n^1(\mathbf{Z}_n)$ and $\Psi_n^2(\mathbf{Z}_n)$ in~\eqref{eq:tightnessleftright} are tight. This is done separately for both of them in the sequel. \smallskip

$\bullet\ \mathbf{\Psi_n^1}$. 
The crucial idea to prove tightness of $\Psi_n^1$ is to define an event with high probability on which we can bound $|\hat{\Psi}_{n,t,h}^{\hat{c}_n}(\mathbf{Z}_n)-\hat{\Psi}_{n,t,h}^{c}(\mathbf{Z}_n) |$ by a quantity independent of $t,h$. For $\delta,\kappa>0$, 
\[ \left| \frac{1}{c^2}-\frac{1}{\hat{c}_n^2}\right|\1_{\{ |\hat{c}_n^2-c^2|\leq \kappa c^2/\sqrt{n}\}} \1_{\{\hat{c}_n^2\geq \delta c^2\}} \leq \frac{\kappa}{\sqrt{n}c^2\delta} \1_{\{ |\hat{c}_n^2-c^2|\leq \kappa c^2/\sqrt{n}\}} \1_{\{\hat{c}_n^2\geq \delta c^2\}},\]
such that
\begin{align*}
    & \P_c\left( \sup_{(t,h) \in \mathcal{T}_n(K)} \big|\hat{\Psi}_{n,t,h}^{\hat{c}_n}(\mathbf{Z}_n) -\hat{\Psi}_{n,t,h}^{c}(\mathbf{Z}_n)\big|>   C\right) 
    \leq P_1+P_2 +P_3,
\end{align*}
where
\begin{align*}
   P_1& := \P_c \left( |\hat{c}_n^2-c^2| > \kappa c^2/\sqrt{n} \right),\\
   P_2&:=  \P_c \left( \hat{c}_n^2< \delta c^2 \right),\\
   P_3&:= \P_c\left( \sup_{(t,h) \in \mathcal{T}_n(K)}  \frac{\kappa}{\sqrt{2n}\delta c^2} \sum_{k=1}^n \frac{\lambda_k(t,h)}{\Lambda_n(t,h)}  n(Z_{k/n}-Z_{(k-1)/n})^2> C\right).
\end{align*}
Let $\epsilon>0$. First, we will choose $\kappa,\delta>0$ in such a way that $P_1,P_2\leq \epsilon/3$. For the respective choices of these parameters, we then show that $C$ may be chosen such that $P_3\leq \epsilon/3$ and tightness of $\Psi_n^1$ follows.\smallskip

\begin{itemize}
    \item[$\bullet\ \mathbf{P_1}$.] By Markov's inequality, independence of increments of $Z$ and $\E_c[(Z_{k/n}-Z_{(k-1)/n})^2]=c^2/n$,
    \begin{align*}
        \P_c \left( |\hat{c}_n^2-c^2| > \frac{\kappa c^2}{\sqrt{n}} \right) &\leq \frac{n}{\kappa^2 c^4}\E_c\left[\left(\hat{c}_n^2-c^2\right)^2\right] = \frac{1}{\kappa^2}\frac1n \sum_{k=1}^n \E_c\left[ \left( \frac{n}{c^2}(Z_{k/n}-Z_{(k-1)/n})^2 -1\right)^2\right] = \frac{2}{\kappa^2},
    \end{align*} 
    where the last equality used $\textrm{Var}(Y)=2$ for $Y\sim\chi_1^2$. For $\kappa\to\infty$, this term clearly vanishes.\smallskip

    \item[$\bullet\ \mathbf{P_2}$.] Let $Y_n\sim\chi_n^2$. For $a>0$ we have the  Chernoff-type bound
    \[ \P\left(Y_n\leq n\delta\right) \leq e^{a n\delta}\E\left[ e^{-a Y_n}\right] = e^{a n\delta} (1+2a)^{-n/2},\]
    where the last step uses that the moment generating function $M_{Y_n}(t)$ of $Y_n$ is $M_{Y_n}(t)=(1-2t)^{-n/2}$ for $t<1/2$ (see p.64 in~\cite{Davison}). For $0<\delta<1$ and $a=(1-\delta)/(2\delta)$ this yields
    \[ \P\left(Y_n\leq n\delta\right) \leq \exp\left(-\frac{n}{2}\left(\delta-1-\log(\delta)\right)\right).\]
    As $\sum_{k=1}^n n(Z_{k/n}-Z_{(k-1)/n})^2\sim \chi_n^2$ under $\P_1$ we obtain
    \[ \P_c\left( \hat{c}_n^2 < \delta c^2\right) = \P_1\left( \sum_{k=1}^n n(Z_{k/n}-Z_{(k-1)/n})^2 \leq n\delta\right) \leq \exp\left(-\frac12\left( \delta-1 -\log(\delta)\right)\right).\]
    As $\delta-1-\log(\delta)\rightarrow\infty$ as $\delta\searrow 0$, we can choose $\delta>0$ small enough such that this probability bound is smaller than $\epsilon/3$.

    \item[$\bullet\ \mathbf{P_3}$.] Here, we start with a deterministic bound using Cauchy-Schwarz' inequality for sums. Recalling $\Lambda_n(t,h)^2=\sum_{k=1}^n \lambda_k(t,h)^2$, we find
    \begin{align*}
        &\sup_{(t,h)\in\mathcal{T}_n(K)}\frac{\kappa}{\sqrt{2n}\delta c^2} \frac{1}{\Lambda_n(t,h)} \sum_{k=1}^n \lambda_k(t,h) n(Z_{k/n}-Z_{(k-1)/n})^2\\
        &\hspace{1cm}\leq \sup_{(t,h)\in\mathcal{T}_n(K)}\frac{\kappa}{\sqrt{2n}\delta c^2} \frac{1}{\Lambda_n(t,h)} \sqrt{\Lambda_n(t,h)^2} \sqrt{\sum_{k=1}^n n^2(Z_{k/n}-Z_{(k-1)/n})^4}\\
        &\hspace{1cm} = \frac{\kappa}{\sqrt{2n}\delta c^2} \sqrt{\sum_{k=1}^n n^2(Z_{k/n}-Z_{(k-1)/n})^4}.
    \end{align*}
    Hence, Markov's inequality and $\E_c[(Z_{k/n}-Z_{(k-1)/n})^4] =3c^4/n^2$ reveal
    \begin{align*}
        P_3 &\leq \P_c\left( \frac{\kappa}{\sqrt{2n}\delta c^2} \sqrt{\sum_{k=1}^n n^2(Z_{k/n}-Z_{(k-1)/n})^4} >C\right) \\
        &\leq \frac{\kappa^2}{2n\delta^2 c^4C^2} \E_c\left[\sum_{k=1}^n n^2(Z_{k/n}-Z_{(k-1)/n})^4 \right] = \frac{3\kappa^2}{2\delta^2 C^2} \stackrel{C\to\infty}{\longrightarrow} 0.
    \end{align*}
\end{itemize}

$\bullet\ \mathbf{\Psi_n^2}$.
We introduce the process $(Y_n(t,h))_{(t,h)\in\mathcal{T}_n(K)}$ with
\[ Y_n(t,h) = \frac{1}{\sqrt{2}}\sum_{k=1}^n \lambda_k(t,h)\left( c^{-2}n\left( Z_{k/n}-Z_{(k-1)/n}\right)^2-1\right)\]
and a semimetric $d_n$ on $\mathcal{T}_n(K)$ via
\[ d_n((t,h),(t',h')) =  \sqrt{\sum_{k=1}^n \left( \lambda_k(t,h) - \lambda_k(t',h')\right)^2}. \]
It is easily checked by Cauchy-Schwarz' inequality that $|\Lambda_n(t,h)-\Lambda_n(t',h')|\leq d_n((t,h),(t',h'))$ for all $(t,h),(t',h')\in\mathcal{T}_n(K)$. Moreover, we define
\[ G_n'(\eta,\delta) :=  \sqrt{2\eta}+ \frac{\sqrt{2}\eta}{\max\{1,\sqrt{n}\delta\}} \]
for which we observe $\sup_{n\in\N}\sup_{\eta\geq 0, 0<\delta\leq 1} G_n'(\eta,\delta)/(1+\eta)<\infty$. Then, by Theorem~$10$ in~\cite{Rohde2009} with the subsequent remark replaced by
\[ G_n'(2\log(1/\delta)+C\log\log(e/\delta),\delta) \leq G_n(2\log(1/\delta),\delta) + o(1)\quad\textrm{ as }\delta\searrow 0\]
and Theorem~$7$ in~\cite{Dumbgen_report} with its subsequent remark, we obtain that 
\[ \left( \sup_{(t,h)\in\mathcal{T}_n(K)} \left( \big|\hat{\Psi}_{n,t,h}^c(\mathbf{Z}_n)\big| - G_n\left( 2\log(1/\Lambda_n(t,h)), \Lambda_n(t,h)\right)\right) \right)_{n\in\N}\]
is tight, once we can prove the following three assertions:
\begin{itemize}
    \item[(i)] For any $(t,h)\in\mathcal{T}_n(K)$ with $\Lambda_n(t,h)\geq\delta$, we have
    \[ \P_c\left( |Y_n(t,h)| \geq \Lambda_n(t,h) G_n'(\eta,\delta)\right) \leq 2\exp\left(-\eta\right).\]
    \item[(ii)] There exist constants $c_1,c_2>0$ such that for all $(t,h),(t',h')\in\mathcal{T}_n(K)$
    \[ \P_c\left( |Y_n(t,h)-Y_n(t',h')|\geq c_1 d_n((t,h),(t',h'))\eta\right) \leq c_2\exp\left(-\eta\right).\]
    \item[(iii)] There exists a constant $c_3>0$ such that
    \[ D\left( (\delta u), \{(t,h)\in\mathcal{T}_n(K): \Lambda_n(t,h)\leq\delta\},d_n\right) \leq c_3 u^{-4}\delta^{-2},\]
    for all $u,\delta\in(0,1]$, where $D(u,A,d_n)$ denotes the capacity number of the set $A$ w.r.t. the semimetric $d_n$, i.e. $D(u,A,d_n) = \max\{ \# A_0: A_0\subset A, d_n(s,t)>u \textrm{ for different } s,t\in A_0\}$. Observe that this is sufficient for an application of Theorem~$10$ in~\cite{Rohde2009} as $N(u,A,d_n)\leq D(u,A,d_n)\leq N(u/2,A,d_n)$ for the covering number $N(u,A,d_n)$, see p.147 in~\cite{Vaart/Wellner}.
\end{itemize}
In what follows, we are going to prove (i)-(iii), thereby establishing tightness of $\Psi_n^2$ and thus finishing the proof of Theorem~\ref{thm:tightness}.\medskip

    \textbf{Proof of (i).} We apply Lemma~\ref{lemma:Bernstein_chi_squared} with $a_k=\lambda_k(t,h)/\sqrt{2\Lambda_n(t,h)^2}$. For this choice, we obviously have
    \[ V_n = \sum_{k=1}^n a_k^2 = \frac12\sum_{k=1}^n \frac{\lambda_k(t,h)^2}{\Lambda_n(t,h)^2} =\frac12.\]
    As $\Lambda_n(t,h)\geq\delta$ by assumption and $\Lambda_n(t,h)^2\geq 1/n$ by Lemma~\ref{lemma:sigmageq} for $(t,h)\in\mathcal{T}_n(K)$, one obtains the bound
    \[ M_n = \max_{k=1,\dots, n}|a_k| = \max_{k=1,\dots, n} \frac{\sqrt{n}\int_{(k-1)/n}^{k/n} \psi_{t,h}(s)ds}{\sqrt{2} \Lambda_n(t,h)} \leq  \frac{1}{\sqrt{2}\max\{1, \sqrt{n}\delta\}}.\]
    As $(\sqrt{n}(Z_{k/n}-Z_{(k-1)/n})/c)_k$ is an i.i.d. sequence of standard Gaussian random variables, we can apply Lemma~\ref{lemma:Bernstein_chi_squared} with $V_n=1/2$ and the bound $M_n\leq \max\{\sqrt{2},\sqrt{2n}\delta\}^{-1}$ which yields~(i).\smallskip

    \textbf{Proof of (ii).} Writing out the definition of $Y_n(t,h)$, we observe
    \[ Y_n(t,h)-Y_n(t',h') = \frac{1}{\sqrt{2}}\sum_{k=1}^n \left(\lambda_k(t,h)-\lambda_k(t',h')\right)\left( c^{-2}n\left( Z_{k/n}-Z_{(k-1)/n}\right)^2-1\right).\]
    Put $b_k = (\lambda_k(t,h)-\lambda_k(t',h'))/\sqrt{2d_n((t,h),(t',h'))^2}$. Then, we find 
    \[ V_n'=\sum_{k=1}^n b_k^2 =1/2\leq 1 \quad\textrm{ and }\quad M_n'=\max_{k=1,\dots,n}|b_k|\leq \sqrt{V_n'}\leq 1, \]
    such that Lemma~\ref{lemma:Bernstein_chi_squared} yields for $\eta\geq 1$ the tail probability bound
    \[ \P_c\left(\frac{|Y_n(t,h)-Y_n(t',h')|}{d_n((t,h),(t',h'))} > 4\eta\right) \leq \P_c\left(\frac{|Y_n(t,h)-Y_n(t',h')|}{d_n((t,h),(t',h'))} > 2\sqrt{\eta}+ 2\eta\right) \leq 2\exp\left(-\eta\right). \]
    For $0<\eta<1$ we trivially bound the probability by one. Hence~(ii) holds true with $c_1=4$ and $c_2=2e$. \smallskip

    \textbf{Proof of (iii).} In this proof, $\lambda(A)$ denotes the Lebesgue measure for a measurable subset $A\subset\R$ and $TV(\psi)$ the variation of $\psi$ (which by assumption is finite). We partition $[0,1]$ into intervals $(M_k)_{1 \leq k \leq N}$ such that $\lambda(M_k) \leq \frac{1}{2TV(\psi)^2}  (u \delta)^2$ with equality for the first $N - 1$. 
    Then we have
    \begin{align*}
        &1= \lambda[0,1]= \sum_{k=1}^N \lambda(M_k) = (N-1) \frac{1}{2TV(\psi)^2} (u \delta)^2+ \lambda(M_N) \ge \frac{(N-1) (u \delta)^2 }{2 TV(\psi)^2},
    \end{align*}
    which yields the bound $N \leq 1+ 2TV(\psi)^2(u\delta)^{-2}$. Next, take $(t,h) \in \mathcal{T}_n(K)$ such that $\Lambda_n(t, h)^2\le \delta^2$. By Lemma~\ref{lemma:sigmageq}, this implies $Ch \leq \delta^2$ for some constant $C>0$. Then, for $t-h\in M_i$ and $t+h \in M_j$, 
    \begin{align*}
    &(j-i-1)\frac{1}{2TV(\psi)^2} (u\delta)^2 \leq 2h \leq \frac{2\delta^2}{C}  \quad  \Leftrightarrow \quad  (j-i)\leq 1 + \frac{4TV(\psi)^2}{Cu^2}.
    \end{align*}
    This implies that there are at most
    \[ \left(1+ \frac{2TV(\psi)^2}{(u\delta)^2} \right)\left( 2 + \frac{4TV(\psi)^2}{Cu^2}\right) \leq \left(2+4TV(\psi)^2+\frac{4TV(\psi)^2}{C}+\frac{8TV(\psi)^4}{C}\right)u^{-4}\delta^{-2} \]
    such pairs $(i,j)$. Let $\mathcal{T}_n' \subset \mathcal{T}_n(K)$ be a maximal subset of $\{(t, h) \in \mathcal{T}_n (K): \Lambda_n (t, h) \leq \delta\}$ with $d_n((t, h), (t', h')) > u \delta$ for arbitrary $(t, h), (t', h') \in \mathcal{T}_n'$. Then, the proof is finished, if we could show that for all $(i, j)$ there is at most one such point $(t_{ij}, h_{ij})$ in $\mathcal{T}_n'$ with $t_{ij} - h_{ij} \in M_i$ and $t_{ij} + h_{ij} \in M_j$. To this aim, we pick $(t, h), (t', h') \in \mathcal{T}_n(K)$ such that $t - h, t' - h' \in M_i$ and $t + h, t' + h' \in M_j$, and are done if we show $d_{n}((t, h), (t', h'))^2 \leq (u \delta)^2$. By Jensen's inequality,
    \begin{align*}
    \begin{split}
        d_n\left((t,h);(t',h')\right)^2   &= \sum_{k=1}^n n\left(\int_{(k-1)/n}^{k/n} \psi_{t,h}(s) - \psi_{t',h'}(s)ds\right)^2 \leq \int_0^1 \left(\psi_{t,h}(s)-\psi_{t',h'}(s)\right)^2ds.
        \end{split}
    \end{align*}
    As $\psi$ is of finite variation, we can write $\psi(x)=\int_{[-1,x]} gdP $ for all but at most countably many numbers $ x \in [-1,1]$ with $P$ a probability measure on $[-1,1]$ and $g$ a measurable function with $|g| \leq TV(\psi)$. Inserting this representation into the last bound and consecutively applying Jensen's inequality and Fubini's theorem reveals
    \begin{align*}
        d_n\left((t,h);(t',h')\right)^2 &\leq \int_0^1 \left(\int_{-1}^{(s-t)/h} g(r)dP(r)- \int_{-1}^{(s-t')/h'} g(r)dP(r)\right)^2 ds\\
        & \leq TV(\psi)^2  \int_{-1}^1 \int_0^1 \1_{B_s((t,h),(t',h'))}(r) ds dP(r),
    \end{align*}
    where $B_s((t,h),(t',h'))=[\min\{(s-t)/h,(s-t')/h'\},\max\{(s-t)/h,(s-t')/h'\}]$.
    Now, as $r\in B_s((t,h),(t',h'))$ if and only if either $rh'+t'\leq s\leq rh+t$ or $rh+t\leq s\leq rh'+t'$, 
    \[ \int_{0}^1 \1_{B_s((t,h),(t',h'))}(r) ds =\left| t'+rh' - t - rh\right| \leq \lambda([t-h,t+h]\triangle[t'-h',t'+h']),\]
    where $\triangle $ denotes the symmetric difference of two sets and the estimate can be seen as $| t'+rh' - t - rh|$ is maximized for $r=\mathrm{sign}(t-t')\mathrm{sign}(h-h')$. Consequently,
    \[ d_n\left((t,h);(t',h')\right)^2 \leq TV(\psi)^2  \lambda([t-h,t+h]\triangle[t'-h',t'+h']) \leq TV(\psi)^2 2 \lambda(M_k) \leq (u\delta)^2.\] 
\end{proof}

\subsection{Proof of the minimax lower bound Theorem~\ref{thm:lower_bound}}\label{Subsection:proof_lower_bound}
First, we present the proofs of the two auxiliary Lemmas~\ref{lemma:hypotheses_Hölder} and~\ref{lemma:hypotheses_solution}. Afterwards, a proof of Theorem~\ref{thm:lower_bound} is provided.

\begin{proof}[Proof of Lemma~\ref{lemma:hypotheses_Hölder}]
    Recall $\psi_\beta\in\mathcal{H}(\beta,1)$. Let $g_{n,j}=\sigma_{n,j}-1$, then for $x,y\in [0,1]$,
    \[ \left| g_{n,j}(x) - g_{n,j}(y)\right| = \frac{L}{2}(1-\epsilon_n) (h_n^{b_n})^\beta \left| \psi_\beta\left( \frac{x-t_{n,j}}{h_n^{b_n}}\right)- \psi_\beta\left( \frac{y-t_{n,j}}{h_n^{b_n}}\right)\right|  \leq \frac{L}{2} (1-\epsilon_n) |x-y|^\beta\]
    and as $|\psi_\beta(x)|\leq 1$ for all $x\in\R$,
    \[ \left| g_{n,j}(x)^2- g_{n,j}(y)^2\right| \leq L(1-\epsilon_n)(h_n^{b_n})^\beta\left| g_{n,j}(x) - g_{n,j}(y)\right| \leq \frac{L^2}{2}(1-\epsilon_n)^2 (h_n^{b_n})^\beta |x-y|^\beta.  \]
    Hence, we obtain $2g_{n,j}\in\mathcal{H}(\beta,L(1-\epsilon_n))$ and $g_{n,j}^2\in\mathcal{H}(\beta,L^2(1-\epsilon_n)^2(h_n^{b_n})^\beta/2)$. This further entails
    \[ \frac{\sigma_{n,j}^2}{\|\sigma_{n,j}\|_{L^2}^2} = \frac{1+2g_{n,j} + g_{n,j}^2}{\|\sigma_{n,j}\|_{L^2}^2} \in \mathcal{H}\left(\beta, \frac{L(1-\epsilon_n)( 1+L(1-\epsilon_n)(h_n^{b_n})^\beta/2)}{\|\sigma_{n,j}\|_{L^2}^2}\right) \subset \mathcal{H}(\beta,L),\]
    where the inclusion is true for $n$ sufficiently large for the following reason: Since $\psi_\beta\geq 0$ for $0<\beta\leq 1$ we have $\|\sigma_{n,j}\|_{L^2}^2\geq 1$ and for $n$ sufficiently large,
     \[ L(1-\epsilon_n)\left( 1+L(1-\epsilon_n)(h_n^{b_n})^\beta/2 \right) \leq L- L\epsilon_n \left( 1- L(h_n^{b_n})^\beta\epsilon_n^{-1}\right) \leq L.\]
    Here, the second inequality follows since $(h_n^{b_n})^\beta\sqrt{\log(n)}\rightarrow 0$ and $\epsilon_n\sqrt{\log(n)}\rightarrow \infty$ as $n\to\infty$ by our choice of $\epsilon_n$.
\end{proof}

\begin{proof}[Proof of Lemma~\ref{lemma:hypotheses_solution}]
    The strategy for constructing $b_n$ is to prove existence of a solution to the equation in~\eqref{eq_proof:hypotheses_solution1} via the implicit function theorem.
    As $\psi_\beta(0)=1$, we have $\sigma_{n,j}^b(t_{n,j}) = 1+L(1-\epsilon_n)(h_n^b)^\beta/2$ and, moreover,
    \[ \|\sigma_{n,j}^b\|_{L^2}^2 = 1 + L(1-\epsilon_n) (h_n^b)^{\beta+1}\int_\R \psi_\beta(x)dx + \frac{L^2}{4}(1-\epsilon_n)^2 (h_n^b)^{2\beta+1} \|\psi_\beta\|_{L^2}^2.\]
    Inserting all of this and using that $(h_n^b)^\beta = (h_n^1)^\beta b^{-\beta/(2\beta+1)} = c_*\rho_n b^{-\beta/(2\beta+1)}/L$ gives
    \[ \left|\frac{\sigma_{n,j}^b(t_{n,j})^2}{\|\sigma_{n,j}^b\|_{L^2}^2} -1\right| = \frac{c_*(1-\epsilon_n)\rho_n b^{-\beta/(2\beta+1)}}{\|\sigma_{n,j}^b\|_{L^2}^2} \left( 1- h_n^b\int_\R \psi_\beta(x)dx + \frac{L}{4}(1-\epsilon_n)(h_n^b)^{\beta}\left(1-h_n^b\|\psi_\beta\|_{L^2}^2\right)\right)\]
    such that the equality in~\eqref{eq_proof:hypotheses_solution1} is equivalent to 
    \[ b^{-\beta/(2\beta+1)}\left( 1- h_n^b\int_\R \psi_\beta(x)dx + \frac{L}{4}(1-\epsilon_n)(h_n^b)^{\beta}\left( 1-h_n^b\|\psi_\beta\|_{L^2}^2\right)\right) - \|\sigma_{n,j}^b\|_{L^2}^2 =0.\]
    By standard algebraic manipulations we obtain that this equality is equivalent to
    \[ F\left((h_n^1)^\beta,\epsilon_n (h_n^1)^\beta, b \right) =0,\]
    where for $\varphi_\beta(u)=\mathrm{sign}(u) |u|^{1/\beta}\in\mathcal{C}^1(\R)$
\begin{align*}
F(u,y,z)&= z^{-\frac{\beta}{2\beta+1}}-1 +\Bigl(-\varphi_\beta(u)\,z^{-\frac{\beta+1}{2\beta+1}}+L\,z^{-\frac{\beta+1}{2\beta+1}}\bigl(y-u\bigr)\varphi_\beta(u)\Bigr)\int_{\mathbb R}\psi_\beta(x)\,dx\\
&\hspace{1cm} +\frac L4\Bigl(1-\varphi_\beta(u)\,z^{-\frac{1}{2\beta+1}}\|\psi_\beta\|_{L^2}^2\Bigr)\Bigl(u\,z^{-\frac{2\beta}{2\beta+1}}-y\,z^{-\frac{2\beta}{2\beta+1}}\Bigr) -\frac{L^2}{4}\|\psi_\beta\|_{L^2}^2\,\bigl(u-y\bigr)^2\,\varphi_\beta(u)\,z^{-1}.
\end{align*}
For this function, we have $F(0,0,1)=0$ and $\frac{\partial}{\partial z}F(0,0,1)=-\beta/(2\beta+1)\neq 0$. By the implicit function theorem, we conclude that there exists an open neighborhood $U\times V\subset\R^2\times \R$ of $(0,0,1)\in\R^3$ such that $F(x,y,f(x,y))=0$ for some continuously differentiable function $f:U\rightarrow V$ with $f(0,0)=1$. This still holds true on any compact subset $\overline{U}\times \overline{V}\subset U\times V$ containing $(0,0,1)$ on which $f$ is Lipschitz continuous as the continuous gradient $\nabla f$ is bounded on this domain. As $(h_n^1)^\beta\rightarrow 0$ and $\epsilon_n (h_n^1)^\beta\rightarrow 0$, we can take $n\geq n_0$ large enough such that $((h_n^1)^\beta, \epsilon_n(h_n^1)^\beta)\in \overline{U}$ and define a solution $b=b_n$ of~\eqref{eq_proof:hypotheses_solution1} as $b_n=f((h_n^1)^\beta,\epsilon_n (h_n^1)^\beta)$. Finally, denote by $C_L$ the Lipschitz constant of the restricted function $f|_{\overline{U}}$ w.r.t. to the Euclidean norm on $\R^2$. Then
\begin{align*}
\log(n) \left| b_n-1\right| = \log(n)\left| f\left( (h_n^1)^\beta, \epsilon_n(h_n^1)^\beta\right) - f\left(0,0\right)\right| \leq C_L\log(n)(h_n^1)^\beta\sqrt{ 1 + \epsilon_n^2} \longrightarrow 0. 
\end{align*}
\end{proof}

\begin{proof}[Proof of Theorem~\ref{thm:lower_bound}]
As $\phi_n$ is a test that uniformly over $H_0$ is of level $\alpha$, we obtain for hypotheses $\sigma_{n,j}$ satisfying $\sigma_{n,j}^2/\|\sigma_{n,j}\|_{L^2}^2\in\mathcal{H}(\beta,L)$ and $d_J(\sigma_{n,j},H_0) \geq (1-\epsilon_n)c_*\rho_n$, 
\begin{align*}
\inf_{ \substack{\sigma^2/\|\sigma\|_{L^2}^2\in\mathcal{H}(\beta,L):\\ d_J(\sigma, H_0)\geq (1-\epsilon_n)c_*\rho_n}} \E_\sigma[\phi_n(\mathbf{Z}_n)] - \alpha 
& \leq \min_{j=1,\dots, N_n} \E_{\sigma_{n,j}}[\phi_n(\mathbf{Z}_n)] - \sup_{\sigma\in H_0} \E_\sigma[\phi_n(\mathbf{Z}_n)] \\
& \leq \frac{1}{N_n} \sum_{j=1}^{N_n} \E_{\sigma_{n,j}}[\phi_n(\mathbf{Z}_n)] - \E_1[\phi_n(\mathbf{Z}_n)]  \leq \E_1 \left[ \left| \frac{1}{N_n} \sum_{j=1}^{N_n} \frac{\textrm{d}\P_{\sigma_{n,j}}}{\textrm{d}\P_1}(\mathbf{Z}_n)  -1\right|\right].
\end{align*}
In what follows, we will choose $N_n$ and $\sigma_{n,j}$ in such a way that the right-hand side converges to zero. To this aim, we are going to construct $\sigma_{n,j}$ such that $\textrm{d}\P_{\sigma_{n,i}}/\textrm{d}\P_1$ and $\textrm{d}\P_{\sigma_{n,j}}/\textrm{d}\P_1$ are independent for $i\neq j$. Then, we obtain from Proposition~$8.1$ in~\cite{Brutsche/Rohde} that for any $\epsilon>0$ and $\gamma\in (0,1]$,
\begin{align*}
    \E_1 \left[ \left| \frac{1}{N_n} \sum_{j=1}^{N_n} \frac{\textrm{d}\P_{\sigma_{n,j}}}{\textrm{d}\P_1}(\mathbf{Z}_n)  -1\right|\right] \leq \sqrt{\epsilon} + 2\epsilon^{-\gamma} N_n^{-(1+\gamma)} \sum_{j=1}^{N_n}\E_1 \left[\left(\frac{\textrm{d}\P_{\sigma_{n,j}}}{\textrm{d}\P_1}(\mathbf{Z}_n)\right)^{1+\gamma}\right]
\end{align*}   
and the claim of Theorem~\ref{thm:lower_bound} follows once we prove convergence of the second summand to zero. To this aim, we are going to prove
\begin{align}\label{eq_proof:lower3}
    \max_{j=1,\dots, N_n}\log\left( N_n^{-\gamma}\E_1 \left[\left(\frac{\textrm{d}\P_{\sigma_{n,j}}}{\textrm{d}\P_1}(\mathbf{Z}_n)\right)^{1+\gamma}\right] \right) \stackrel{n\to\infty}{\longrightarrow} -\infty.
\end{align}
To prove this, we recall the representation~\eqref{eq:likelihood_ratio} of the log-likelihood ratio, $(n(Z_{k/n}-Z_{(k-1)/n})^2)_{k=1,\dots,n}$ being a vector of independent $\chi_1^2$-distributed random variables, and the moment generating function $M_Y(t)=(1-2t)^{-1/2}$ for $t<1/2$ for $Y\sim\chi_1^2$ (see p.64 in~\cite{Davison}). If $\|\sigma_{n,j}\|_{\sup}\rightarrow 1$ (at rate independent of $j$), we can choose $n$ large enough such that $2(1-\|\sigma_{n,j}^2\|_{\sup}^{-1})<1$ and obtain
\begin{align*}
    \E_1 \left[\left(\frac{\textrm{d}\P_{\sigma_{n,j}}}{\textrm{d}\P_1}(\mathbf{Z}_n)\right)^{1+\gamma}\right] &= \prod_{k=1}^n \left( n\int_{(k-1)/n}^{k/n} \sigma_{n,j}(s)^2 ds\right)^{-(1+\gamma)/2}  \\
    &\hspace{0.5cm}\cdot\E_1\left[\exp\left( \frac12 (1+\gamma) n\left(Z_{(k-1)/n}-Z_{k/n}\right)^2 \left( 1- \frac{1}{n\int_{(k-1)/n}^{k/n} \sigma_{n,j}(s)^2 ds}\right)\right) \right] \\
    &= \prod_{k=1}^n \left( n\int_{(k-1)/n}^{k/n} \sigma_{n,j}(s)^2 ds\right)^{-\frac{1+\gamma}{2}} \left( 1-(1+\gamma)\left( 1- \frac{1}{n\int_{(k-1)/n}^{k/n} \sigma_{n,j}(s)^2 ds}\right)\right)^{-\frac12}
\end{align*}
such that after some elementary calculations in the first step and Jensen's inequality applied to the convex function $x\mapsto -\log(x)$ in the second,
\begin{align}\label{eq_proof:lower2}
\begin{split}
    &\log\left(\E_1 \left[\left((\textrm{d}\P_{\sigma_{n,j}}/\textrm{d}\P_1)(\mathbf{Z}_n)\right)^{1+\gamma}\right] \right)\\
    &\hspace{0.5cm} = -\frac{1}{2}\sum_{k=1}^n \left[\gamma\log\left( n\int_{(k-1)/n}^{k/n} \sigma_{n,j}(s)^2ds\right) + \log\left(-\gamma n\int_{(k-1)/n}^{k/n}\sigma_{n,j}^2(s)ds + 1+\gamma\right) \right] \\
    &\hspace{0.5cm} \leq -\frac12\sum_{k=1}^n\left[ \gamma n\int_{(k-1)/n}^{k/n} \log\left(\sigma_{n,j}(s)^2\right) ds + n\int_{(k-1)/n}^{k/n} \log\left( -\gamma\sigma_{n,j}(s)^2+1+\gamma\right) ds \right] \\
    &\hspace{0.5cm} = \frac{n}{2} \int_0^1 \left( -\gamma\log\left(\sigma_{n,j}(s)^2\right) - \log\left( 1+\gamma-\gamma\sigma_{n,j}(s)^2\right) \right) ds.
\end{split}
\end{align}
From now on, we have to distinguish the cases $\beta\leq 1$, for which we even obtain the optimal constant and $\beta>1$ for which only rate optimality is shown. As outlined after Theorem~\ref{thm:lower_bound}, for obtaining the optimal constant, the hypotheses $\sigma_{n,j}$ have to be chosen more carefully. The remaining structure of the proof is as follows: We first construct suitable hypotheses for $\beta\leq 1$ and then prove~\eqref{eq_proof:lower3} for these, where the derivation continues from the last line of~\eqref{eq_proof:lower2}. This will finish the proof in case $\beta\leq 1$. Afterwards, we provide the construction of the hypotheses in case $\beta>1$ and again establish~\eqref{eq_proof:lower3} based on~\eqref{eq_proof:lower2}. \smallskip

\textbf{Construction of hypotheses for $\beta\leq 1$.} We give the details for the construction of the hypotheses for $\beta \leq 1$ which was already outlined in Section~\ref{Section:power_properties} right after the statement of Theorem~\ref{thm:lower_bound}. Let $J=[J_1,J_2]$ and for a parameter $b>0$ we define $h_n^b$ as in~\eqref{eq:h_nb} and $\sigma_{n,j}^b$ as in~\eqref{eq:sigma_j} such that it remains to specify $t_{n,j}^b\in [J_1,J_2]$. Define \vspace{-0.1cm}
\[ N_n^b:=\left\lfloor  \frac{J_2-J_1}{2(h_n^b+1/n)}\right\rfloor\quad\textrm{ and }\quad  t_{n,j}^b := J_1+ (2j-1)\left( h_n^b + \frac1n\right) \quad \textrm{ for } j=1,\dots, N_n^b,\]
such that in particular $J_1+h_n^b\leq t_{n,j}^b\leq J_2-h_n^b$ for all $j=1,\dots, N_n^b$. By Lemma~\ref{lemma:hypotheses_Hölder} there exists $n_0$ such that for all $n\geq n_0$ we have $(\sigma_{n,j}^b)^2/\|\sigma_{n,j}^b\|_{L^2}^2\in\mathcal{H}(\beta,L)$ for all $1\leq j\leq N_n^b$. Moreover, by Lemma~\ref{lemma:hypotheses_solution} there exists a sequence $(b_n)_{n\in\N}$ with $\log(n)|b_n-1|\rightarrow 0$ such that $\sigma_{n,j}^{b_n}$ satisfies $d_J(\sigma_{n,j}^{b_n},H_0)\geq (1-\epsilon_n)c_*\rho_n$ for $n\geq n_1$. In the sequel, we always assume $n\geq \max\{n_0,n_1\}$ and set for ease of notation
\[ t_{n,j}:=t_{n,j}^{b_n}, \qquad N_n:=N_n^{b_n}\quad \textrm{ and } \quad \sigma_{n,j}:=\sigma_{n,j}^{b_n}.\] 
It is clear by definition~\eqref{eq:sigma_j} of $\sigma_{n,j}$ that $\max_{1\leq j\leq N_n} \|\sigma_{n,j}-1\|_{\sup}\rightarrow 0$ as $n\to\infty$. Moreover, the distance between two consecutive $t_{n,j}$ is $2h_n^{b_n}+2/n$, while the support of $\psi_\beta((s-t_{n,j})/h_n^{b_n})$ has length $2h_n^{b_n}$. Hence, there does not exist an index $k$ for which $(k-1)/n$ is contained in the support of $\psi_\beta((s-t_{n,i})/h_n^{b_n})$ and $k/n$ in that of $\psi_\beta((s-t_{n,j})/h_n^{b_n})$ for $i<j$. Because the increments of $Z$ are independent, this reveals that $\textrm{d}\P_{\sigma_{n,i}}/\textrm{d}\P_1$ and  $\textrm{d}\P_{\sigma_{n,j}}/\textrm{d}\P_1$ are independent for $i\neq j$.\smallskip

\textbf{Proof of~\eqref{eq_proof:lower3} for $\beta\leq 1$.}
Set $g_{n,j}=\sigma_{n,j}-1$. Then $\max_j\|g_{n,j}\|_{\sup}\rightarrow 0$ as $n\to\infty$ and for $n$ large enough such that $b_n\geq 1/2 $,
\begin{align}\label{eq_proof:lower1}
\begin{split}
    n\int_0^1 g_{n,j}(s)^m ds 
    &\leq L^m2^{-m}\|\psi_\beta\|_{\sup}^m \left( \frac{c_*}{L}\right)^{\frac{m\beta+1}{\beta}} \left( 2\log(n)\right)^{\frac{m\beta+1}{2\beta+1}}  n^{-\frac{(m-2)\beta}
    {2\beta+1}},
\end{split}
\end{align} 
which converges to zero for $m\geq 3$. In the next step, we will continue from~\eqref{eq_proof:lower2} by using the inequality $-\log(1+x)\leq -x+x^2/2$ valid for $x\geq 0$ and a third order Taylor expansion in $x=0$ with Lagrange remainder that yields the equality $-\log(1-x)=x+x^2/2+x^3/(3(1-\xi)^3)$ for some $\xi$ between $0$ and $x$. We apply the Taylor expansion to $x=\gamma(g_{n,j}(s)^2 + 2g_{n,j}(s))$ where the intermediate value $\xi_{n,j}(s)$ lies between $0$ and $\gamma(g_{n,j}(s)^2 + 2g_{n,j}(s))$. In particular, we can choose $n$ large enough such that $1/(1-\xi_{n,j})^3 \leq 2$ uniformly in $s$ as $\|g_{n,j}\|_{\sup}\rightarrow 0$. With~\eqref{eq_proof:lower1}, we obtain the following estimate for the right-hand side of~\eqref{eq_proof:lower2}:
\begin{align*}
  & n\int_0^1 \left( -\gamma\log\left(\sigma_{n,j}(s)\right) - \frac12\log\left( 1+\gamma-\gamma\sigma_{n,j}(s)^2\right) \right) ds \\
    &\hspace{0.5cm} \leq n\int_0^1 \left[-\gamma g_{n,j}(s) + \frac{\gamma}{2} g_{n,j}(s)^2 + \frac12\left( \gamma g_{n,j}(s)^2 + 2\gamma g_{n,j}(s) + \frac12 \gamma^2 \left( g_{n,j}(s)^2 + 2g_{n,j}(s)\right)^2 \right.\right.\\
    &\hspace{2.5cm} \left.\left. + \frac{\gamma^3 \left( g_{n,j}(s)^2 + 2g_{n,j}(s)\right)^3}{3(1-\xi_{n,j}(s))^3}\right)\right] ds  \\
    &\hspace{0.5cm} = n \gamma (1+\gamma) \int_0^1 g_{n,j}(s)^2 ds + o(1). 
\end{align*}
Inserting this result into~\eqref{eq_proof:lower2} and using the definition of $g_{n,j}$ and the fact that $N_n\geq (J_2-J_1)/(4h_n^{b_n})$ for $n$ large enough then yields 
\begin{align*}
    &\log\left( N_n^{-\gamma}\E_1 \left[\left((\textrm{d}\P_{\sigma_{n,j}}/\textrm{d}\P_1)(\mathbf{Z}_n)\right)^{1+\gamma}\right] \right)\\
    &\hspace{0.5cm} \leq -\gamma \log\left( N_n\right) + n \gamma (1+\gamma) \int_0^1 g_{n,j}(s)^2 ds + o(1) \\
    &\hspace{0.5cm} \leq -\gamma\log\left(\frac{1}{h_n^{b_n}}\right) + \frac14 L^2 \gamma(1+\gamma)(1-\epsilon_n)^2 n(h_n^{b_n})^{2\beta} \int_0^1 \psi_\beta\left( \frac{s-t_{n,j}}{h_n^{b_n}}\right)^2 ds + \mathcal{O}(1)\\
    &\hspace{0.5cm} =\frac{\gamma}{2\beta+1}\left( \log(\log(n)) - \log(n)\right) + \frac14 L^2 \gamma(1+\gamma)(1-\epsilon_n)^2 \left(\frac{c_*}{L}\right)^{\frac{2\beta+1}{\beta}} \frac{\log(n)}{b_n} \|\psi_\beta\|_{L^2}^2 +\mathcal{O}(1) \\
    &\hspace{0.5cm} = \frac{\gamma}{2\beta+1}\left( \log(\log(n)) - \log(n)\right) + \gamma(1+\gamma)(1-\epsilon_n)^2 \frac{1}{2\beta+1} \frac{\log(n)}{b_n} + \mathcal{O}(1),
\end{align*}
where the last step uses the definition of the optimal constant $c_*$ in~\eqref{eq:optimal_constant}. Finally, we set $\gamma=\epsilon_n$. Then for all $n$ large enough such that the last display is valid, we obtain \newpage
\begin{align*}
    \log\left( N_n^{-\gamma}\E_1 \left[\left(\frac{\textrm{d}\P_{\sigma_{n,j}}}{\textrm{d}\P_1}(\mathbf{Z}_n)\right)^{1+\gamma}\right] \right) &= - \frac{\epsilon_n^2\left( 1+ \mathcal{O}(\epsilon_n)\right)}{2\beta+1}\log(n) + \frac{\epsilon_n}{2\beta+1}\log(\log(n))\\
    &\hspace{1cm} + \mathcal{O}(\epsilon_n)\log(n)\left(\frac{1}{b_n}-1\right) + \mathcal{O}(1). 
\end{align*}
Note that this bound does not depend on $j$. Recalling $\log(n)|b_n-1|\rightarrow 0$ by Lemma~\ref{lemma:hypotheses_solution}, the third summand vanishes. Finally, as
\[ -\frac{\epsilon_n^2(1+\mathcal{O}(\epsilon_n))}{2\beta+1}\log(n) + \frac{\epsilon_n}{2\beta+1}\log(\log(n)) = -\frac{\epsilon_n^2\log(n)}{2\beta+1} \left( 1 +\mathcal{O}(\epsilon_n) - \frac{\log(\log(n))}{\epsilon_n\log(n)}\right) \rightarrow -\infty\]
by our assumption $\epsilon_n\sqrt{\log(n)}\rightarrow\infty$, the convergence in~\eqref{eq_proof:lower3} actually holds true and Theorem~\ref{thm:lower_bound} is proven for $\beta\leq 1$.\smallskip

\textbf{Construction of hypotheses for $\beta> 1$.}
Let $J=[J_1,J_2]$. We denote the support of $\psi_\beta$ by $[-R_\beta,R_\beta]$ and consider hypotheses of the form \vspace{-0.1cm}
\[ \sigma_{n,j}^2=1+L(1-\epsilon_n)h_n^{\beta}\psi_{\beta}\left(\frac{s-t_{n,j}}{h_n}\right),\]
where
\begin{align}\label{eq:alternativebetage1}
     h_n=\left(\frac{4^{\left(\beta/(2\beta+1)\right)}c}{L}\right)^{\frac{1}{\beta}}\left(\frac{\log(n)}{n}\right)^{1/(2\beta+1)} \quad\textrm{ and }\quad c=\left(\frac{L^{1/\beta}}{(2\beta+1)\|\psi_{\beta}\|^2_{L^2}}\right)^{\beta/(2\beta+1)}.
\end{align}
Their number $N_n$ and the location parameter $t_{n,j}$ are defined as
\[ N_n = \left\lfloor \frac{J_2-J_1}{2R_\beta(h_n+1/n)}\right\rfloor \quad\textrm{ and }\quad t_{n,j} := J_1+(2j-1)R_\beta\left( h_n + \frac1n\right) \quad \textrm{ for } j=1,\dots, N_n.\]
Analogously to the case $\beta\leq 1$ we have $J_1+R_\beta h_n\leq t_{n,j}\leq J_2- R_\beta h_n$ for all $j=1,\dots, N_n$, that the likelihood ratios $\textrm{d}\P_{\sigma_{n,i}}/\textrm{d}\P_1$ and $\textrm{d}\P_{\sigma_{n,j}}/\textrm{d}\P_1$ are independent for $i\neq j$ and $\max_{j=1,\dots, N_n}\|\sigma_{n,j}^2-1\|_{\sup}\rightarrow 0$.\\
Clearly, $\sigma_{n,j}^2$ is $\lfloor\beta\rfloor$-times differentiable as $\psi_\beta$ is and we observe
\begin{align*}
    \left|(\sigma_{n,j}^2)^{(\lfloor \beta \rfloor)}(x)- (\sigma_{n,j}^2)^{(\lfloor \beta \rfloor)}(y)\right|
    &= L(1-\epsilon_n)h_n^{\beta-\lfloor \beta \rfloor} \left| \psi_{\beta}^{(\lfloor \beta \rfloor)}\left(\frac{x-t_{n,j}}{h_n}\right)- \psi_{\beta}^{(\lfloor \beta \rfloor)}\left(\frac{y-t_{n,j}}{h_n}\right)\right|\\
    &\leq L(1-\epsilon_n)|x-y|^{\beta -\lfloor \beta \rfloor},
\end{align*}
such that $\sigma_{n,j}^2\in\mathcal{H}(\beta,L)$. As $ \|\sigma_{n,j}\|_{L^2}^2 = 1 + L(1-\epsilon_n) (h_n)^{\beta+1}\int_\R \psi_\beta(x)dx \geq 1$ because by Theorem~$1$ in~\cite{Leonov} $\int_\R\psi_\beta(x)dx>0$, we conclude $\sigma^2_{n,j}/\| \sigma_{n,j}\|_{L^2}^2 \in \mathcal{H}(\beta,L)$. \\
Finally, we check that $d_J(\sigma_{n,j},H_0)\geq (1-\epsilon_n)c\rho_n$. Here, we first note that $\psi_\beta(0)=1$ such that $\sigma_{n,j}^2(t_{n,j}) = 1+L(1-\epsilon_n)(h_n)^\beta$. With the explicit form of $\|\sigma_{n,j}\|_{L^2}^2$ and $(h_n)^\beta  = 4^{\beta/(2\beta+1)}c\rho_n/L$, we obtain
\[ d_J(\sigma_{n,j},H_0)\geq \left|\frac{\sigma_{n,j}(t_{n,j})^2}{\|\sigma_{n,j}\|_{L^2}^2} -1\right| = \frac{(1-\epsilon_n)c\rho_n}{\|\sigma_{n,j}\|_{L^2}^2} 4^{\beta/(2\beta+1)} \left( 1- h_n\int_\R \psi_\beta(x)dx \right).\]
The right-hand side needs to be greater than $(1-\epsilon_n)c\rho_n$, which is true for $n$ large enough since $h_n\to 0$ and $\max_j(\|\sigma_{n,j}\|_{L^2}^2-1)\rightarrow 0$.\\

\textbf{Proof of~\eqref{eq_proof:lower3} for $\beta> 1$.}
Inequality~\eqref{eq_proof:lower1} remains valid for $\sigma_{n,j}^2-1$ with different constants such that we obtain by $-\log(1+x)\leq -x + x^2/2 + \mathcal{O}(|x|^3)$ for $|x|\leq 1/2$ and~\eqref{eq_proof:lower1} along the lines for $\beta\leq 1$, that
\[ \log\left( \E_1 \left[\left(\frac{\textrm{d}\P_{\sigma_{n,j}}}{\textrm{d}\P_1}(\mathbf{Z}_n)\right)^{1+\gamma}\right] \right) \leq \frac{n}{4} \gamma(1+\gamma) \int_0^1 L^2(1-\epsilon_n)^2h_n^{2\beta}\psi_{\beta}\left(\frac{s-t_{n,j}}{h_n}\right)^2ds + o(1).\]
By definition of the constant $c$ and $h_n$ in~\eqref{eq:alternativebetage1} and the estimate $N_n\geq (J_2-J_1)/(4R_\beta h_n)$ for $n$ large enough, this estimate reveals
\begin{align*}
    &\log\left( N_n^{-\gamma}\E_1 \left[\left(\frac{\textrm{d}\P_{\sigma_{n,j}}}{\textrm{d}\P_1}(\mathbf{Z}_n)\right)^{1+\gamma}\right] \right)\\
     &\hspace{1cm}\leq -\gamma\log\left(\frac{1}{h_n}\right) + \frac14 L^2 \gamma(1+\gamma)(1-\epsilon_n)^2 n(h_n)^{2\beta} \int_0^1 \psi_\beta\left( \frac{s-t_{n,j}}{h_n}\right)^2 ds + \mathcal{O}(1)\\
    &\hspace{1cm} = \frac{\gamma}{2\beta+1}\left( \log(\log(n)) - \log(n)\right) + \gamma(1+\gamma)(1-\epsilon_n)^2 \frac{1}{2\beta+1} \log(n) + \mathcal{O}(1)
\end{align*}
and~\eqref{eq_proof:lower3} follows by the same arguments as it did for $\beta\leq 1$ with the choice $\gamma=\epsilon_n$.
\end{proof}

\subsection{Proof of minimax upper bounds}\label{Subsection:Proof_upper_bound}
This Subsection contains the proofs of Theorem~\ref{thm:upper_bound} and Theorem~\ref{thm:adaptivity} in the respective order. Before we start with these, we state one preliminary lemma together with its proof.

\begin{lemma}\label{lemma:Hölder_density_sup_bounded}
    Let $\beta,L>0$. Then there exists a constant $C=C(\beta,L)>0$, which is bounded on compact
    subsets of $(0,\infty)^2$, such that for every function $f:[0,1]\rightarrow\R$ with $f^2/\|f\|_{L^2}^2\in\mathcal{H}(\beta,L)$ we have
    \[ \|f\|_{\sup}^2 \leq C\|f\|_{L^2}^2.\]
\end{lemma}
\begin{proof}
Set $g=f^2/\|f\|_{L^2}^2$, so that $g\geq 0$, $\int_0^1 g(x)dx=1$ and $g\in\mathcal{H}(\beta,L)$. Let $m=\lfloor\beta\rfloor$, fix $x_0\in[0,1]$ and let $P$ denote the Taylor polynomial of $g$ of degree $m$ at $x_0$. By the integral form of the Taylor remainder and the H\"older continuity of $g^{(m)}$, $R:=g-P$ satisfies 
\[ |R(x)|\leq \frac{L}{m!}|x-x_0|^\beta\leq \frac{L}{m!},\qquad x\in[0,1], \]
where for $m=0$ this simply reads $|g(x)-g(x_0)|\leq L$. Hence $\|P\|_{L^1[0,1]}\leq \|g\|_{L^1[0,1]}+L/m! = 1+L/m!$. As the space $\Pi_m$ of polynomials of degree at most $m$ is finite-dimensional and $\|\cdot\|_{L^1[0,1]}$ is a norm on $\Pi_m$, there exists $M_m<\infty$ with $\|p\|_{\sup}\leq M_m\|p\|_{L^1[0,1]}$ for all $p\in\Pi_m$. Consequently,
\[ \|g\|_{\sup}\leq \|P\|_{\sup}+\|R\|_{\sup}\leq M_m\Big(1+\frac{L}{m!}\Big)+\frac{L}{m!} =:C(\beta,L). \]
Since $m$ takes only finitely many values for $\beta$ in a compact set and $C(\beta,L)$ is non-decreasing in $L$, the constant $C$ is bounded on compact subsets of $(0,\infty)^2$.
\end{proof}

\begin{proof}[Proof of Theorem~\ref{thm:upper_bound}]
For any $\sigma:[0,1]\rightarrow (0,\infty)$ and $k=1,\dots, n$ we define
\[ I_k(\sigma) = n\int_{(k-1)/n}^{k/n} \sigma(s)^2 ds \qquad\textrm{ and }\qquad \hat{c}_{n,\sigma}^2 = \sum_{k=1}^n I_k(\sigma)(Z_{k/n}-Z_{(k-1)/n})^2.\]
The first idea in the proof is to apply a change of measure and subsequently work under $\P_1$. By the construction of the test $T_n$ we find for any $t\in J$ and $h>0$,
\begin{align*}
    \P_\sigma\left( T_n(\mathbf{Z}_n) =1\right) &\geq \P_\sigma\left( \big|\hat{\Psi}_{n,t,h}^{\hat{c}_n}\big| > \kappa_n^\alpha + G_n\left(2\log(1/\Lambda_n(t,h)), \Lambda_n(t,h)\right) \right) \\
    &= \P_1\Bigg( \left| \frac{1}{\sqrt{2}\Lambda_n(t,h) \hat{c}_{n,\sigma}^2}\sum_{k=1}^n \lambda_k(t,h) \left(I_k(\sigma)n(Z_{k/n}-Z_{(k-1)/n})^2 - \hat{c}_{n,\sigma}^2\right) \right|\\
    &\hspace{5cm}  > \kappa_n^\alpha + G_n\left(2\log(1/\Lambda_n(t,h)), \Lambda_n(t,h)\right)\Bigg),
\end{align*} 
where the second step uses the definition of $\hat{\Psi}_{n,t,h}^{\hat{c}_n}$ and $\hat{c}_n^2$ given in~\eqref{eq:local_statistic_c} and~\eqref{eq:c_hat} together with the fact that the vector of increments $((Z_{k/n}-Z_{(k-1)/n})^2)_k$ is distributed under $\P_\sigma$ as $(I_k(\sigma)(Z_{k/n}-Z_{(k-1)/n})^2)_k$ is under $\P_1$. It remains to evaluate this probability, where we note that $\E_1[\hat{c}_{n,\sigma}^2] = \|\sigma\|_{L^2}^2$ and $\E_1[n(Z_{k/n}-Z_{(k-1)/n})^2]=1$ such that by adding and subtracting expectations, 
\begin{align}\label{eq_proof:upper1}
    \P_\sigma\left( T_n(\mathbf{Z}_n) =1\right) \geq \P_1\Big( |A_n^\sigma(t,h)| > \kappa_n^\alpha + |B_n^\sigma(t,h)| + G_n\left(2\log(1/\Lambda_n(t,h)), \Lambda_n(t,h)\right) \Big),
\end{align} 
where we used $|a+b|\geq |a|-|b|$ and abbreviated
\begin{align*}
    A_n^\sigma(t,h) &= \frac{\frac{1}{\sqrt{2}}\sum_{k=1}^n \lambda_k(t,h)\left( I_k(\sigma) - \|\sigma\|_{L^2}^2\right)}{\Lambda_n(t,h) \|\sigma\|_{L^2}^2} ,\\
    B_n^\sigma(t,h) &= \frac{\frac{1}{\sqrt{2}}\sum_{k=1}^n \lambda_k(t,h)\left( I_k(\sigma) - \|\sigma\|_{L^2}^2\right)}{\Lambda_n(t,h)}\left(\frac{1}{\hat{c}_{n,\sigma}^2} - \frac{1}{\|\sigma\|_{L^2}^2}\right) \\
    &\hspace{0.5cm} + \frac{\frac{1}{\sqrt{2}}\sum_{k=1}^n \lambda_k(t,h) \left(I_k(\sigma)(n(Z_{k/n}-Z_{(k-1)/n})^2-1)+ \|\sigma\|_{L^2}^2 - \hat{c}_{n,\sigma}^2\right)}{\Lambda_n(t,h) \hat{c}_{n,\sigma}^2}. 
\end{align*}
We now want to show that the right-hand side of~\eqref{eq_proof:upper1} converges to one for any sequence $(\sigma_n)_{n\in\N}$ with $\sigma_n^2/\|\sigma_n\|_{L^2}^2\in\mathcal{H}(\beta,L)$ and $d_J(\sigma_n,H_0)\geq (1+\epsilon_n)c\rho_n$ and appropriate $(t,h)=(t_n,h_n)$. Let $(\sigma_n)_{n\in\N}$ be such a sequence. The location $t_n\in J$ and the bandwidth $h_n$,
which will be a multiple of $\rho_n^{1/\beta}$, are specified below, as their choice differs
for $\beta\leq 1$ and $\beta>1$.
In order to prove convergence of the right-hand side of~\eqref{eq_proof:upper1} towards one for $(t_n,h_n)$, we use that $(\kappa_n^\alpha)_{n\in\N}$ is bounded by Theorem~\ref{thm:tightness}. Then, this convergence follows once we establish
\begin{itemize}
    \item[(1)] $|A_n^{\sigma_n}(t_n,h_n)|-G_n(2\log(1/\Lambda_n(t_n,h_n)),\Lambda_n(t_n,h_n)) \rightarrow \infty$, and
    \item[(2)] the sequence $(|B_n^{\sigma_n}(t_n,h_n)|)_{n\in\N}$  is tight under $\P_1$.
\end{itemize}
Note that the quantity in (1) is deterministic. It will turn out that this term requires distinguishing between $\beta\leq 1$ and $\beta>1$, whereas one can argue the same way for all $\beta>0$ in (2). Subsequently, we prove items (1) and (2) which finishes the proof of Theorem~\ref{thm:upper_bound} by~\eqref{eq_proof:upper1}. \smallskip

\textbf{Item (1) for $\beta\leq 1$.} Here, we use weights $\lambda_k(t,h)$ based on the optimal recovery kernel $\psi_\beta$ and define
\begin{align}\label{eq:h_detect}
    h_n = \left(\frac{(1+\epsilon_n)c_*\rho_n}{L}\right)^{\frac{1}{\beta}}= \left(\frac{(1+\epsilon_n)c_*}{L}\right)^{\frac{1}{\beta}} \left( \frac{\log(n)}{n}\right)^{\frac{1}{2\beta+1}}.
\end{align} 
Set $f_n=\sigma_n^2/\|\sigma_n\|_{L^2}^2-1$. By assumption on $\sigma_n$, we have $f_n\in\mathcal{H}(\beta,L)$ and, as $d_J(\sigma_n,H_0)\geq (1+\epsilon_n)c_*\rho_n$ and $f_n$ is
continuous, we can choose $t_n\in J$ with $|f_n(t_n)|\geq (1+\epsilon_n)c_*\rho_n$. We assume w.l.o.g.\ that $f_n(t_n)>0$, otherwise replace $f_n$ by $-f_n$. Then, by the H\"older assumption,
for any $s\in [t_n-h_n, t_n+h_n]$ it holds that
\begin{align*}
    f_n(s) \geq (1+\epsilon_n)c_*\rho_n-L|s-t_n|^\beta = (1+\epsilon_n)c_*\rho_n \left( 1-\frac{|s-t_n|^\beta}{h_n^\beta}\right) = (1+\epsilon_n)c_*\rho_n\psi_\beta\left(\frac{s-t_n}{h_n}\right),
\end{align*}
where the last step is valid due to the representation~\eqref{eq:optimal_recovery_kernel} of $\psi_\beta$. Moreover, the second last step implies the rough bound $f_n(s)\geq -c_*\rho_n$ for $s\in [t_n-h_n-1/n,t_n+h_n+1/n]$. Using the first property, we obtain for every $k$ with $t_n-h_n\leq (k-1)/n \leq k/n\leq t_n+h_n$,
\[ \frac{I_k(\sigma_n) -\|\sigma_n\|_{L^2}^2}{\|\sigma_n\|_{L^2}^2} = n\int_{(k-1)/n}^{k/n} \left( \frac{\sigma_n(s)^2}{\|\sigma_n\|_{L^2}^2} -1\right) ds \geq (1+\epsilon_n)c_*\rho_n n\int_{(k-1)/n}^{k/n} \psi_\beta\left(\frac{s-t_n}{h_n}\right) ds.\]
Together with the definition of $\lambda_k(t_n,h_n)$ for $\psi_\beta$ and the definition of $\Lambda_n(t_n,h_n)$ we obtain
\begin{align*}
    \big|A_n^{\sigma_n}(t_n,h_n)\big| &\geq \frac{1}{\sqrt{2}}\sum_{k=\lceil n(t_n-h_n)\rceil+1}^{\lfloor n(t_n+h_n)\rfloor} \frac{\lambda_k(t_n,h_n)}{\Lambda_n(t_n,h_n)} (1+\epsilon_n)c_*\rho_n n\int_{(k-1)/n}^{k/n} \psi_\beta\left( \frac{s-t_n}{h_n}\right) ds - \sqrt{\frac{2}{n}}\frac{c_*\rho_n}{\Lambda_n(t_n,h_n)} \\
    &= \frac{(1+\epsilon_n)c_*\rho_n \sqrt{n}}{\sqrt{2}} \sum_{k=\lceil n(t_n-h_n)\rceil +1}^{\lfloor n(t_n+h_n)\rfloor} \frac{\lambda_k(t_n,h_n)^2}{\Lambda_n(t_n,h_n)}- \sqrt{\frac{2}{n}}\frac{c_*\rho_n}{\Lambda_n(t_n,h_n)} \\
    &\geq \frac{(1+\epsilon_n)c_*\rho_n \sqrt{n}}{\sqrt{2}}\Lambda_n(t_n,h_n) - \frac{\sqrt{2}(2+\epsilon_n)c_*\rho_n}{\sqrt{n}\Lambda_n(t_n,h_n)},
\end{align*}
where we added the boundary summands $k=\lceil n(t_n-h_n)\rceil, \lfloor n(t_n+h_n)\rfloor +1$ in the last step (for all other $k$ we have $\lambda_k(t_n,h_n)=0$) and subtracted an upper bound on them using $\lambda_k(t_n,h_n)^2\leq 1/n$. By Lemma~\ref{lemma:sigmageq}, we know $\Lambda_n(t_n,h_n)^2\geq Ch_n$ and because $nh_n\to\infty$,
\begin{align}\label{eq_proof:upper8}
    \big|A_n^{\sigma_n}(t_n,h_n)\big| &\geq \frac{(1+\epsilon_n)c_*\rho_n \sqrt{nh_n}}{\sqrt{2}}\|\psi_\beta\|_{L^2} -\frac{1+\epsilon_n}{\sqrt{2}}c_*\rho_n\sqrt{n}\left| \Lambda_n(t_n,h_n) - \sqrt{h_n}\|\psi_\beta\|_{L^2}\right| + o(1).
\end{align}
By substitution, $h_n\|\psi_\beta\|_{L^2}^2$ = $\|(\psi_\beta)_{t_n,h_n}\|_{L^2}^2$ and by the explicit representation for $0<\beta\leq 1$ one reads off $\|\psi_\beta\|_{\sup}\leq 1$. By the intermediate value theorem, we find variables $(k-1)/n\leq \eta_k,\tilde{\eta}_k\leq k/n$ such that
\begin{align*}
    \left| \Lambda_n(t_n,h_n)^2- h_n\|\psi_\beta\|_{L^2}^2\right| &\leq \sum_{k=1}^n \left| n\left(\int_{(k-1)/n}^{k/n} (\psi_\beta)_{t_n,h_n}(s) ds\right)^2 - \int_{(k-1)/n}^{k/n} (\psi_{\beta})_{t_n,h_n}(s)^2 ds\right| \\
    &= \frac1n \sum_{k=1}^n \left| (\psi_{\beta})_{t_n,h_n}(\eta_k)^2-(\psi_{\beta})_{t_n,h_n}(\tilde{\eta}_k)^2 \right|\\
    &\leq \frac2n \sum_{k=1}^n  \left| (\psi_{\beta})_{t_n,h_n}(\eta_k)-(\psi_{\beta})_{t_n,h_n}(\tilde{\eta}_k) \right|\\
    &\leq \frac{2}{(nh_n)^\beta}\frac1n \sum_{k=0}^n \1_{[t_n-h_n, t_n+h_n]}(k/n) \\
    &\leq \frac{2}{(nh_n)^\beta}\frac1n \left( 2nh_n + 1\right) \leq \frac{6h_n}{(nh_n)^\beta},
\end{align*}
where the last inequality is valid for $n$ large enough such that $2nh_n+1\leq 3nh_n$. Because for $a,b\geq 0$ we have $|\sqrt{a}-\sqrt{b}| \leq |a-b|/\sqrt{b}$, we obtain for these $n$,
\begin{align*}
    2c_*\rho_n\sqrt{n}\left| \Lambda_n(t_n,h_n) - \sqrt{h_n}\|\psi_\beta\|_{L^2}\right| &\leq 12L  (1+\epsilon_n)^{-1}\|\psi_\beta\|_{L^2}^{-1} \sqrt{nh_n}n^{-\beta} \\
    &= 12L (1+\epsilon_n)^{-1}\|\psi_\beta\|_{L^2}^{-1}\left(\frac{(1+\epsilon_n)c_*}{L}\right)^{\frac{1}{2\beta}} \log(n)^{\frac{1}{2(2\beta+1)}} n^{-\frac{2\beta^2}{2\beta+1}}
\end{align*} 
which converges to zero, and hence~\eqref{eq_proof:upper8} yields
\begin{align}\label{eq_proof:upper2}
    \big|A_n^{\sigma_n}(t_n,h_n)\big| \geq \frac{(1+\epsilon_n)c_*\rho_n \sqrt{nh_n}}{\sqrt{2}}\|\psi_\beta\|_{L^2} + o(1).
\end{align} 
Next, we turn to the term that involves $G_n$. As $\Lambda_n(t_n,h_n)^2\geq Ch_n$ by Lemma~\ref{lemma:sigmageq},
\begin{align}\label{eq_proof:lower_correction_negligible}
    \frac{\log(1/\Lambda_n(t_n,h_n))}{\sqrt{n} \Lambda_n(t_n,h_n)} \leq \frac{\log(1/h_n)+\log(1/C)/2}{\sqrt{Cn h_n}} \rightarrow 0,
\end{align} 
as $\log(1/h_n) = \mathcal{O}(\log(n))$ and $nh_n/\log(n)^2\rightarrow\infty$. With the definition of $h_n$ in~\eqref{eq:h_detect}, Lemma~\ref{lemma:sigmageq} and~\eqref{eq_proof:lower_correction_negligible}, it follows that
\begin{align}\label{eq_proof:upper4}
\begin{split}
    G_n(2\log(1/\Lambda_n(t_n,h_n)),\Lambda_n(t_n,h_n)) &= \sqrt{4\log(1/\Lambda_n(t_n,h_n))} + o(1) \\
    &= \sqrt{-2\log(h_n)+\mathcal{O}(1)} + o(1) \\
    &= \sqrt{\frac{2}{2\beta+1}\log(n)} + o(1).
\end{split}
\end{align}
Combining~\eqref{eq_proof:upper2} and~\eqref{eq_proof:upper4} and inserting the definition of $\rho_n$, $h_n$ and $c_*$ finally reveals
\begin{align*}
    &\big|A_n^{\sigma_n}(t_n,h_n)\big|-G_n(2\log(1/\Lambda_n(t_n,h_n)),\Lambda_n(t_n,h_n))   \\
    &\hspace{1cm} \geq \frac{(1+\epsilon_n)c_*\rho_n \sqrt{nh_n}}{\sqrt{2}}\|\psi_\beta\|_{L^2} - \sqrt{\frac{2}{2\beta+1}\log(n)}+ \mathcal{O}(1) \\
    &\hspace{1cm} = \frac{\|\psi_\beta\|_{L^2}}{\sqrt{2}L^{1/(2\beta)}} \sqrt{n} (1+\epsilon_n)^{\frac{2\beta+1}{2\beta}} \left(\frac{\log(n)}{n}\right)^\frac12 \left(\frac{4L^{1/\beta}}{(2\beta+1)\|\psi_\beta\|_{L^2}^2}\right)^\frac12 - \sqrt{\frac{2}{2\beta+1}\log(n)}+ \mathcal{O}(1) \\
    &\hspace{1cm} \geq \epsilon_n\sqrt{\log(n)} \sqrt{\frac{2}{2\beta+1}} + \mathcal{O}(1)
\end{align*}
where in the last step we applied the bound $(1+\epsilon_n)^{(2\beta+1)/(2\beta)}\geq 1+\epsilon_n$. The resulting term now clearly converges to $\infty$ by our assumption on $\epsilon_n$.\smallskip

\textbf{Item (1) for $\beta> 1$.} Set $f_n=\sigma_n^2/\|\sigma_n\|_{L^2}^2-1$. Then $f_n\in\mathcal{H}(\beta,L)$ and, by Lemma~\ref{lemma:Hölder_density_sup_bounded}, $\|f_n\|_{\sup}\leq D:=C(\beta,L)+1$. Let $x_n\in J$ be a maximizer of $|f_n|$ over $J$, so that $\delta_n:=|f_n(x_n)|=d_J(\sigma_n,H_0)\geq(1+\epsilon_n)c\rho_n$. Applying Lemma~\ref{lem:rohde2} to $f_n/D\in\mathcal{H}(\beta,L/D)$ on $J$ yields a constant $c'=c'(\beta,L)>0$ such that $|f_n|\geq\delta_n/2$ on $J\cap[x_n-c'\delta_n^{1/\beta},x_n+c'\delta_n^{1/\beta}]$. Define 
\[ h_n:=\frac{c'}{2}\big((1+\epsilon_n)c\rho_n\big)^{1/\beta}, \]
so that $2h_n\leq c'\delta_n^{1/\beta}$. For $n$ large enough such that $4h_n\leq J_2-J_1$, the interval $J$ contains $[x_n,x_n+2h_n]$ or $[x_n-2h_n,x_n]$, and both are contained in $[x_n-c'\delta_n^{1/\beta},x_n+c'\delta_n^{1/\beta}]$. Let $t_n:=x_n+h_n$ in the first case and $t_n:=x_n-h_n$ in the second one. As $2h_n\leq c'\delta_n^{1/\beta}$, we have $[t_n-h_n,t_n+h_n]\subset J$ and $|f_n|\geq\delta_n/2\geq(1+\epsilon_n)c\rho_n/2$ on this interval. By continuity, $f_n$ has a constant sign there and we may assume w.l.o.g.\ $f_n>0$, see the argument for $\beta\leq 1$. We now choose weights $\lambda_k(t,h)$ built on a kernel $\psi$ satisfying the assumptions of Theorem~\ref{thm:tightness}. Hence, for any $t_n-h_n\leq (k-1)/n\leq k/n\leq t_n+h_n$,
\begin{align*}
    \frac{I_k(\sigma_n) - \|\sigma_n\|_{L^2}^2}{\|\sigma_n\|_{L^2}^2} = n\int_{(k-1)/n}^{k/n}
    f_n(s)\, ds \geq \frac{(1+\epsilon_n)c\rho_n}{2}.
\end{align*}
Moreover, by Jensen's inequality,
\[ \Lambda_n(t_n,h_n)^2 = \sum_{k=1}^n n\left(\int_{(k-1)/n}^{k/n} \psi_{t_n,h_n}(s)ds\right)^2 \leq \sum_{k=1}^n \int_{(k-1)/n}^{k/n} \psi_{t_n,h_n}(s)^2ds = h_n \|\psi\|_{L^2}^2\]
and we obtain (following the steps to derive~\eqref{eq_proof:upper2})
\begin{align*}
    \big| A_n^{\sigma_n}(t_n,h_n)\big| &\geq \frac{1}{\sqrt{8}}\sum_{k=\lceil n(t_n-h_n)\rceil+1}^{\lfloor n(t_n+h_n)\rfloor} \frac{(1+\epsilon_n)c\rho_n }{\Lambda_n(t_n,h_n)} \sqrt{n}\int_{(k-1)/n}^{k/n} \psi_{t_n,h_n}(s)ds + \mathcal{O}(1)\\
    &\geq (1+\epsilon_n)\rho_n \sqrt{n h_n} \frac{c}{\sqrt{8}\|\psi\|_{L^2}}\int_\R \psi(s)ds + \mathcal{O}(1).
\end{align*}
The observation~\eqref{eq_proof:upper4} remains unchanged and thus
\begin{align*}
    &\big|A_n^{\sigma_n}(t_n,h_n)\big|-G_n(2\log(1/\Lambda_n(t_n,h_n)),\Lambda_n(t_n,h_n))   \\
    &\hspace{1cm} \geq (1+\epsilon_n)\rho_n \sqrt{n h_n} \frac{c}{\sqrt{8}\|\psi\|_{L^2}}\int_\R \psi(s)ds -\sqrt{\frac{2}{2\beta+1}\log(n)}+ \mathcal{O}(1) \\
    &\hspace{1cm} = (1+\epsilon_n)^{\frac{2\beta+1}{2\beta}} \sqrt{\log(n)} c^{\frac{2\beta+1}{2\beta}} \frac{(c'/2)^{1/2}}{\sqrt{8}\|\psi\|_{L^2}}\int_\R\psi(s)ds -\sqrt{\frac{2}{2\beta+1}\log(n)}+ \mathcal{O}(1).
\end{align*}
This term goes to $\infty$ if we choose the constant $c$ large enough such that
\[ c^{\frac{2\beta+1}{2\beta}} \frac{(c'/2)^{1/2}}{\sqrt{8}\|\psi\|_{L^2}}\int_\R\psi(s)ds \geq \sqrt{\frac{2}{2\beta+1}}.\]

\medskip
\textbf{Item (2).} We first bound
\[ \left| B_n^{\sigma_n}(t_n,h_n)\right| \leq T_n^1 + T_n^2 + T_n^3 + T_n^4,\]
and prove tightness of each $T_n^j$ for $j=1,2,3,4$, where the four terms are given as
\begin{align*}
    T_n^1 &:= B_n^{\sigma_n}(t_n,h_n)\1_{\{ \hat{c}_{n,\sigma_n}^2 < \|\sigma_n\|_{L^2}^2/2\}}, \\
    T_n^2 &:= \frac{\sqrt{2}\sum_{k=1}^n \lambda_k(t_n,h_n)| I_k(\sigma_n) - \|\sigma_n\|_{L^2}^2 |}{\Lambda_n(t_n,h_n)}\frac{|\hat{c}_{n,\sigma_n}^2-\|\sigma_n\|_{L^2}^2|}{\|\sigma_n\|_{L^2}^4},\\
    T_n^3 &:= \frac{\sqrt{2}\left|\sum_{k=1}^n \lambda_k(t_n,h_n) I_k(\sigma_n)(n(Z_{k/n}-Z_{(k-1)/n})^2-1)\right|}{\Lambda_n(t_n,h_n) \|\sigma_n\|_{L^2}^2}, \\
    T_n^4 &:= \frac{\sqrt{2}\sum_{k=1}^n \lambda_k(t_n,h_n) \big| \|\sigma_n\|_{L^2}^2 - \hat{c}_{n,\sigma_n}^2 \big|}{\Lambda_n(t_n,h_n) \|\sigma_n\|_{L^2}^2}.
\end{align*}
Let $\epsilon>0$. Subsequently, we prove $\P_1(|T_n^j|>K)\leq \epsilon$ for $j=1,2,3,4$ and $K>0$ large enough separately. In all cases, we use the bound
\begin{align}\label{eq_proof:upper6}
\begin{split}
    \E_1\left[ \left( \|\sigma_{n}\|_{L^2}^2 - \hat{c}_{n,\sigma_n}^2\right)^2\right] &= \frac{1}{n^2} \sum_{k=1}^n I_k(\sigma_n)^2 \E_1\left[ \left(n(Z_{k/n}-Z_{(k-1)/n})^2-1\right)^2\right]  \\
    &= 2\sum_{k=1}^n \left(\int_{(k-1)/n}^{k/n} \sigma_n(s)^2 ds\right)^2\leq \frac2n\|\sigma_n\|_{\sup}^2 \|\sigma_n\|_{L^2}^2,
\end{split}
\end{align}
where the first equality is a consequence of the independence of increments of the process $Z$ and the second one uses $\mathrm{Var}(Y)=2$ for $Y\sim\chi_1^2$.
\begin{itemize}
    \item[$\bullet\ \mathbf{T_n^1.}$] By Markov's inequality and equation~\eqref{eq_proof:upper6},
    \[ \P_1\left(|T_n^1|>K\right) \leq \P_1\left( \hat{c}_{n,\sigma_n}^2 < \|\sigma_n\|_{L^2}^2/2\right) \leq \P_1\left(\left| \|\sigma_{n}\|_{L^2}^2 - \hat{c}_{n,\sigma_n}^2\right| > \|\sigma_n\|_{L^2}^2/2\right) \leq \frac{8}{n}\frac{\|\sigma_n\|_{\sup}^2}{\|\sigma_n\|_{L^2}^2}, \]
    which by Lemma~\ref{lemma:Hölder_density_sup_bounded} converges to zero. 
\end{itemize}
To prove tightness of $T_n^j$ for $j=2,3,4$ we show that the respective second moment under $\P_1$ is bounded by a polynomial in $\|\sigma_n\|_{\sup}/\|\sigma_n\|_{L^2}$ which by Lemma~\ref{lemma:Hölder_density_sup_bounded} is bounded uniformly in $n$. The claim then follows by an application of Markov's inequality.
\begin{itemize}
    \item[$\bullet\ \mathbf{T_n^2.}$] By definition,
    \[ \frac{I_k(\sigma_n)-\|\sigma_n\|_{L^2}^2}{\|\sigma_n\|_{L^2}^2} = n\int_{(k-1)/n}^{k/n} \left( \frac{\sigma_n(s)^2}{\|\sigma_n\|_{L^2}^2} -1\right) ds \leq \frac{\|\sigma_n\|_{\sup}^2}{\|\sigma_n\|_{L^2}^2} +1.\]
    Additionally, we obtain using Jensen's inequality
    \begin{align}\label{eq_proof:upper7}
        \left(\sum_{k=1}^n \lambda_k(t_n,h_n)\right)^2 = n^2 \left(\frac1n \sum_{k=1}^n \lambda_k(t_n,h_n)\right)^2 \leq n\sum_{k=1}^n \lambda_k(t_n,h_n)^2 = n\Lambda_n(t_n,h_n)^2.
    \end{align}
    Inserting all of these inequalities together with the moment bound~\eqref{eq_proof:upper6} yields
    \[ \E_1\left[ |T_n^2|^2\right] \leq 4 \left( \frac{\|\sigma_n\|_{\sup}^2}{\|\sigma_n\|_{L^2}^2} +1\right)^2 \frac{\|\sigma_n\|_{\sup}^2}{\|\sigma_n\|_{L^2}^2}. \]
     \smallskip

    \item[$\bullet\ \mathbf{T_n^3.}$] Here, we use $I_k(\sigma_n) \leq \|\sigma_n\|_{\sup}^2$ and $\E_1[(n(Z_{k/n}-Z_{(k-1)/n})^2 -1)^2] =2$ such that by independence of the increments of $Z$,
    \begin{align*}
        \E_1\left[ |T_n^3|^2\right] &\leq \frac{2}{\Lambda_n(t_n,h_n)^2 \|\sigma_n\|_{L^2}^4} \E_1\left[ \left( \sum_{k=1}^n \lambda_k(t_n,h_n) I_k(\sigma_n)\left( n(Z_{k/n}-Z_{(k-1)/n})^2-1\right)\right)^2\right]\\
        & \leq \frac{4 \|\sigma_n\|_{\sup}^4}{\Lambda_n(t_n,h_n)^2\|\sigma_n\|_{L^2}^4} \sum_{k=1}^n \lambda_k(t_n,h_n)^2 = \frac{4\|\sigma_n\|_{\sup}^4}{\|\sigma_n\|_{L^2}^4}.
    \end{align*} 

    \item[$\bullet\ \mathbf{T_n^4.}$] By the moment bound~\eqref{eq_proof:upper6} and~\eqref{eq_proof:upper7} we directly obtain $\E_1\left[ |T_n^4|^2\right] \leq 4\|\sigma_n\|_{\sup}^2/\|\sigma_n\|_{L^2}^2$.
\end{itemize} 

\end{proof}

\begin{proof}[Proof of Theorem~\ref{thm:adaptivity}]
    The proof follows the lines of that of Theorem~\ref{thm:upper_bound}. Here, the reasoning for Item~(2) remains completely the same by noting that the constant in Lemma~\ref{lemma:Hölder_density_sup_bounded} can be chosen uniformly over $(\beta,L)\in[\beta_1,\beta_2]\times[L_1,L_2]$. To achieve rate adaptivity in both parameters $\beta$ and $L$, we can no longer choose the optimal recovery kernel $\psi_\beta$ for building the test statistic. This means that for $\beta\leq 1$, we have to argue as in the proof of Item~(1) for $\beta>1$. Here, it is important to note that $c'=c'(\beta,L)$ is bounded away from zero over $(\beta,L)\in [\beta_1,\beta_2]\times [L_1,L_2]$ by Lemma~\ref{lem:rohde2}. The same applies for the remainder in~\eqref{eq_proof:upper4} that is uniformly $o(1)$ over this set. With these two ingredients one can follow the proof of Item~(1) for $\beta>1$ step-by-step.\\
    For sharp adaptivity in $L$ for $\beta\leq 1$, one has to note that $\log(c_*/L)$ is uniformly bounded in $L\in [L_1,L_2]$. Then~\eqref{eq_proof:upper4} remains valid and the proof of Item~(1) for $\beta\leq 1$ remains the same as in the proof of Theorem~\ref{thm:upper_bound}.
\end{proof}

\textbf{Funding.} 
L. Riepl gratefully acknowledges support by the DFG Research Unit $5381$, RO $3766/8$-$1$ and CRC $1597$, Project-ID $499552394$.

\bibliographystyle{plainnat}
\bibliography{bibliography}

\end{document}